\documentclass[10pt]{amsart}
\usepackage[british]{babel}
\usepackage{amssymb,mathrsfs,amsfonts,amsmath,amstext,amsopn,mathdesign,bbm,mathtools}
\usepackage{texdraw,graphicx,extarrows,subfigure}%
\usepackage{multirow,hhline,color,colortbl}
\usepackage{mhchem,gensymb}
\usepackage[justification=justified]{caption}
\usepackage{pgfplots,geometry}
\usepackage[all]{xy}
\usetikzlibrary{positioning}

\usepackage[normalem]{ulem}

\usepackage{epstopdf}

\usepackage[numbers,sort&compress]{natbib}
\bibpunct[, ]{[}{]}{,}{n}{,}{,}

\theoremstyle{plain}
\newtheorem{theorem}{Theorem}[section]
\newtheorem{lemma}[theorem]{Lemma}
\newtheorem{corollary}[theorem]{Corollary}
\newtheorem{proposition}[theorem]{Proposition}

\theoremstyle{definition}

\newtheorem{example}[theorem]{Example}

\theoremstyle{remark}
\newtheorem{remark}{Remark}

\newcommand{\eb}{\begin{example}}
\newcommand{\ee}{\end{example}}
\newcommand{\eqb}{\begin{equation}}
\newcommand{\eqe}{\end{equation}}
\newcommand{\spb}{\begin{split}}
\newcommand{\spe}{\end{split}}
\newcommand{\cab}{\begin{cases}}
\newcommand{\cae}{\end{cases}}
\newcommand{\thmb}{\begin{theorem}}
\newcommand{\thme}{\end{theorem}}
\newcommand{\qb}{\begin{question}}
\newcommand{\qe}{\end{question}}
\newcommand{\leb}{\begin{lemma}}
\newcommand{\lee}{\end{lemma}}
\newcommand{\cob}{\begin{corollary}}
\newcommand{\coe}{\end{corollary}}
\newcommand{\enb}{\begin{enumerate}}
\newcommand{\ene}{\end{enumerate}}
\newcommand{\mb}{\begin{matrix}}
\newcommand{\me}{\end{matrix}}
\newcommand{\cenb}{\begin{center}}
\newcommand{\cene}{\end{center}}

\newcommand{\prb}{\begin{proof}}
\newcommand{\pre}{\end{proof}}
\newcommand{\prob}{\begin{proposition}}
\newcommand{\proe}{\end{proposition}}
\newcommand{\rb}{\begin{remark}}
\newcommand{\re}{\end{remark}}

\newcommand{\tS}{\text{S}}
\newcommand{\ka}{\kappa}
\newcommand{\lf}{\lfloor}
\newcommand{\rf}{\rfloor}
\newcommand{\lc}{\lceil}
\newcommand{\rc}{\rceil}

\newcommand{\lt}{\left}
\newcommand{\rt}{\right}

\newcommand{\om}{\omega}

\newcommand{\N}{\mathbb{N}}

\newcommand{\R}{\mathbb{R}}
\newcommand{\Z}{\mathbb{Z}}

\newcommand{\bz}{\mathbf{0}}

\newcommand{\cF}{\mathcal{F}}

\newcommand{\cI}{\mathcal{I}}

\newcommand{\cO}{\mathcal{O}}

\newcommand{\cT}{\Omega}

\newcommand{\cX}{\mathcal{Y}}

\newcommand{\sE}{{\sf E}}

\newcommand{\si}{{\sf i}}
\newcommand{\sK}{{\sf K}}
\newcommand{\sL}{{\sf L}}
\newcommand{\sN}{{\sf N}}

\newcommand{\so}{{\sf o}}
\newcommand{\sP}{{\sf P}}
\newcommand{\sQ}{{\sf Q}}

\newcommand{\sS}{{\sf S}}
\newcommand{\sT}{{\sf T}}

\newcommand{\supp}{{\rm supp\,}}
\newcommand{\spa}{{\rm span\,}}

\newcommand{\wx}{\widetilde{x}}

\newcommand{\ux}{\underline{x}}
\newcommand{\uc}{c_*}
\newcommand{\bc}{c^*}
\newcommand{\uD}{D_*}
\newcommand{\bD}{D^*}
\newcommand{\bu}{u^*}
\newcommand{\srcsize}{\@setfontsize{\srcsize}{5pt}{5pt}}

\makeatletter%
\newcommand\xldownrupharpoon[2][]{%
\ext@arrow 0099{\ldownrupharpoonfill@}{#1}{#2}}
 \def\ldownrupharpoonfill@{%
\arrowfill@\leftharpoondown\relbar\rightharpoonup}
\makeatother

\begin{document}


\title[Structural classification of CMTCs with applications]
{Structural classification of continuous time Markov chains with applications}

\author{Chuang Xu}
\address{
Faculty of Mathematics\\
Technical University of Munich\\
Garching bei M\"{u}nchen, 85748, Germany.}
\email{chuang.xu@ma.tum.de}
\author{Mads Christian Hansen}
\author{Carsten Wiuf}
\address{Department of Mathematical Sciences\\
University of Copenhagen\\ 
2100 K{\o}benhavn, Denmark.}
\email{madschansen@gmail.com}
\email{wiuf@math.ku.dk}

\begin{abstract}
This paper is motivated by examples from \emph{stochastic reaction network theory}. The $Q$-matrix of a stochastic reaction network can be derived from the \emph{reaction graph}, an edge-labelled directed graph encoding the jump vectors of an associated continuous time Markov chain on the invariant space  $\N^d_0$. An open question is how to decompose the   space $\N^d_0$ into neutral, trapping, and escaping states, and open and closed communicating classes, and whether this can be done from the reaction graph alone. Such general continuous time Markov chains can be understood as natural generalizations of birth-death processes, incorporating multiple different birth and death mechanisms. We characterize the structure of $\N^d_0$ imposed by a general  $Q$-matrix   generating continuous time Markov chains with values in $\N^d_0$, in terms of the set of jump vectors and their corresponding transition rate functions. Thus the setting is not limited to stochastic reaction networks. Furthermore, we define structural equivalence of two $Q$-matrices, and provide sufficient conditions for structural equivalence. Examples are abundant in  applications. We apply the results to stochastic reaction networks, a Lotka-Volterra model in ecology, the EnvZ-OmpR system in systems biology, and a class of extended branching processes,  none of which are birth-death processes.
\end{abstract}

\maketitle

\section{Introduction}

Models based on continuous time Markov chains (CTMCs)  are ubiquitous, for example, in genetics \cite{E79}, epidemiology \cite{PCMV15}, ecology \cite{G83}, biochemistry and systems biology \cite{W06}, sociophysics \cite{WH83}, and queueing theory \cite{GH98}.
For CTMCs on a denumerable state space, \emph{structural} classification of the state space is among the fundamental topics and areas of interest.

To investigate the dynamics of a CTMC, we need to know the structure of the state space. If there are no absorbing states, then  the CTMC might be transient, null or positive recurrent, in which case it might admit no stationary measures, a (non-summable) stationary measure, or a non-degenerate ergodic stationary distribution, respectively. Otherwise, if there exist absorbing states, then a degenerate stationary distribution concentrated on the absorbing states trivially exists. In this scenario, the existence of a  quasi-stationary distribution (QSD) may be of interest. Generically, stationary distributions are supported on closed communicating classes while QSDs, if they exist, are concentrated on open communicating classes.

A standard way of defining a class of CTMCs is via a \emph{$Q$-matrix} (the \emph{infinitesimal generator}) on an ambient space $\cX$. Each initial state in $\cX$ defines a CTMC with state space determined by the connectivity of the $Q$-matrix. In applications it is often \emph{difficult} to characterize the structure of the ambient space (its decomposition into communicating classes) purely in terms of the directed graph associated with the $Q$-matrix.  This has for example been revealed in recent work on \emph{stochastic reaction networks} (SRNs)\--- stochastic dynamical models of chemical systems  \cite{ABCJ18,HW20}. SRNs are particular CTMCs on the non-negative integer lattices $\N_0^d$, defined by  $Q$-matrices. In \cite{HW20}, the aim is to identify open and closed communicating classes as means to discuss the existence of stationary and quasi-stationary distributions.  In \cite{ABCJ18}, the interest is  extinction of the process and an algorithmic approach is suggested to identify absorbing states. Thus, a characterization of the absorbing states and the open communicating classes leading to them, is warranted.   In other contexts, it is important to characterize when an SRN is \emph{essential}, that is, when  $\N^d_0$ can be decomposed into a disjoint union of closed communicating classes \cite{PCK14}. However, in general, it is not straightforward  to conclude whether an SRN is essential or not.

SRNs with mass-action kinetics (a particular polynomial form of transition rate functions) can be represented by a (finite) labelled directed graph, called the \emph{reaction graph}, where arrows are reactions and labels are reaction rate constants, see example below.
The reaction graph encodes the possible jump vectors of the CTMC. Despite the simplicity of the graph, it seems non-trivial to deduce the structure of the ambient space $\N_0^d$. This hinges primarily on the fact that all reactions might not be \emph{active} (have non-zero transition rate) on the entire $\N_0^d$, but only on states $y\in\N_0^d$ such that $y\ge x_r$  for some reaction-dependent state $x_r\in\N_0^d$.

Consider the following three SRNs with one species ($\text{S}$) and mass-action kinetics:
\[\begin{tikzpicture}[node distance=2.5em, auto]
 \tikzset{
    pil/.style={
           ->,
           shorten <=2pt,
           shorten >=2pt,}
}
 \node[] (a) {};
 \node[above=1.5em of a] (n0) {\!\!\!\!\!\!A)};
  \node[right=1.2em of a] (n1) {};
  \node[above=-.5em of n1] (m1) {$\ka_1$};
\node[right=of a] (b) {2\tS};
\node[right=1.3em of b] (n2) {};
  \node[above=.9em of b] (m3) {$\ka_3$};
  \node[left=of b] (aa) {\tS} edge[pil, black, bend left=0] (b);
\node[right=of b] (c) {3\tS};
\node[right=of b] (cc) {} edge[pil, black, bend right=40] (aa);
\node[right=.8em of c] (d) {};
 \node[above=1.5em of d] (n00) {\!\!\!\!\!\!B)};
  \node[right=.95em of d] (n11) {};
  \node[above=-.5em of n11] (m11) {$\ka_1$};
\node[right=of d] (e) {2\tS};
\node[right=1em of e] (n22) {};
  \node[above=-.5em of n22] (m22) {$\ka_2$};
  \node[above=.9em of e] (m33) {$\ka_3$};
  \node[left=of e] (dd) {\tS} edge[pil, black, bend left=0] (e);
\node[right=of e] (f) {3\tS};
 \node[left=of f] (ee) {} edge[pil, black, bend left=0] (f);
\node[right=of e] (ff) {} edge[pil, black, bend right=40] (dd);
\node[right=1em of f] (g) {};
 \node[above=1.5em of g] (n000) {\!\!\!\!\!\!C)};
  \node[right=.9em of g] (n111) {};
  \node[above=-.5em of n111] (m111) {$\ka_1$};
\node[right=of g] (h) {2\tS};
\node[right=.8em of h] (n222) {};
  \node[above=-.5em of n222] (m222) {$\ka_2$};
  \node[above=.9em of h] (m333) {$\ka_3$};
  \node[left=of h] (gg) {\tS} edge[pil, black, bend left=0] (h);
\node[right=of h] (i) {3\tS};
  \node[right=1em of i] (n444) {};
  \node[above=-.5em of n444] (m444) {$\ka_4$};
 \node[left=of i] (hh) {} edge[pil, black, bend left=0] (i);
\node[right=of h] (ii) {} edge[pil, black, bend right=40] (gg);
\node[right=of i] (j) {4\tS.};
 \node[left=of j] (ii) {} edge[pil, black, bend left=0] (j);
\end{tikzpicture}\]
The corresponding $Q$-matrices are defined on $\N_0$. A reaction,  $n\text{S}\ce{->[$\kappa$]} m\text{S}$,  encodes  jumps from $x$ to $x+m-n$ with reaction rate function $\lambda(x)=\kappa x(x-1)\ldots(x-n+1)$, $\kappa>0$, such that several reactions might give rise to the same jump vector, for example $\text{S}\ce{->}2\text{S}$ and $2\text{S}\ce{->}3\text{S}$.
The second reaction network is \emph{weakly reversible} (a union of strongly connected components), and is thus essential  \cite{PCK14}, whereas the other two reaction networks are not. However, by the criteria of structural equivalence established in this paper (Theorem~\ref{th-0}), the underlying $Q$-matrices associated with A and C are both \emph{structurally equivalent} to that of B, and hence the three SRNs are all essential.

In this paper, we  study the decomposition of the ambient space $\N_0^d$ into communicating classes, and provide results on the classification of the ambient space into neutral, trapping, escaping and open/closed non-singleton communicating classes, based on a $Q$-matrix defined on $\N^d_0$. The classification is given terms of the set of jump vectors and the support of the transition rate functions, extracted from the $Q$-matrix.  Moreover, we provide necessary and sufficient conditions for the set of absorbing states (neutral and trapping states; typically called the \emph{extinction set} in population process applications) to be non-empty and finite, respectively.

Furthermore, we define the aforementioned structural equivalence for two $Q$-matrices pertaining the decomposition of their ambient spaces, and provide simple checkable conditions for two $Q$-matrices to be structurally equivalent. Such criteria for structural equivalence may provide further insight into identification of reaction networks  with a prescribed structure of the state spaces \emph{via} reaction graphs.

When the CTMCs generated by a $Q$-matrix are restricted to  one-dimensional lattice intervals of $\N^d_0$, we provide complete characterization of all states within such intervals (finite or infinite). In addition, the classification result yields simple checkable criteria for irreducibility of one-dimensional CTMCs. Beyond SRNs, these one-dimensional CTMCs can be understood as \emph{generalized birth-death processes}, with arbitrary jump sizes. For instance, bursty production (eg., of proteins) may induce extra birth mechanisms and the catastrophes may also provide extra death mechanisms (other than natural death). Furthermore, it is found that even within SRNs with mass-action kinetics and $d>1$, there are examples with denumerable irreducible components and very different dynamics on the individual components (for example, explosive \emph{vs.} exponentially ergodic behavior) \cite{XHW20a}. This phenomenon does not appear for $d=1$. Hence, it is necessary to study the delicate structure of the ambient space, in order to fully understand the dynamics of the associated CTMCs. In three companion papers, the classification herein provides means to understand the dynamics of one-dimensional CTMCs \cite{XHW20a,XHW20b} with applications to classification of SRNs \cite{WX20}.

The outline of the paper is as follows: In Section~\ref{sec2}, the notation is introduced and background on CTMCs is reviewed. Main results on structural classification of a $Q$-matrix on  $\N^d_0$ are provided in Section~\ref{sec3}. Applications are in Section~\ref{sec4}. We apply the results to various examples of SRNs, in particular the stochastic Lotka-Volterra system and the biologically important EnvZ-OmpR system, in addition to an extended class of branching processes.

Three lengthy proofs of main results are given in subsequent sections. Finally, some useful elementary propositions and lemmata are appended.

\section{Preliminaries}\label{sec2}

\subsection{Notation}

Let $\R$, $\R_{\ge 0}$, $\R_{>0}$ be the set of real, non-negative real, and positive real  numbers, respectively. Let $\Z$ be the set of integers, $\N=\Z\cap \R_{>0}$ and $\N_0=\N\cup\{0\}$.

For $a\in\R$, let $\lceil a\rceil$ be the ceiling function (i.e., the minimal integer $\ge a$), $\lfloor a\rfloor$ the floor function (i.e., the maximal integer $\le a$), and ${\rm sgn}(a)$ the sign function of $a$. For  $a\in\R\cup\{+\infty\}$, let $\boldsymbol{a}=a(1,\cdots,1)\in\R^d\cup\{\boldsymbol\infty\}$ ($\boldsymbol\infty$ is used if at least one coordinate is $\infty$).
For $x,y\in\R^d$, let $x\ge y$ (similarly, $y\le x$, $x>y$, $y<x$)  if it holds coordinate-wise.

 For $x\in\R^d$, $B\subseteq\R^d$ and $j=1,\cdots,d$, define
$$\widehat x_j=(x_1,\cdots,x_{j-1},x_{j+1},\cdots,x_d),\quad \widehat B_j=\{\widehat x_j\colon x\in B\},\quad B_j=\{x_j\colon x\in B\}.$$

For $a\in\R^d$, $B\subseteq\R$, let $aB=\{ab\colon b\in B\}$. For $b\in\R^d$, $A\subseteq\R^d$,
let $A+b=\{a+b\colon a\in A\}$. Let $\mathbbm{1}_{A}$ denote the indicator function of a set $A$.
For a non-empty subset $A\subseteq\N^d_0$, let
$$\overline{A}=\{x\in \N^d_0\colon x\ge y\,\, \text{for some}\ y\in A\}.$$

The set $A$ is the {\em minimal set} of a non-empty set $B\subseteq\R^d$ if $A$ consists of all elements $x\in B$ such that
$x\not\ge y$ for all $y\in B\setminus\{x\}.$ A set $B$ has a non-empty minimal set if and only if $B$ is bounded from below with respect to the order induced by $\le$. In particular, if $B\subseteq\N^d_0$, then its minimal set $A$ is {\em finite} and $\overline{A}=\overline{B}$ (Proposition~\ref{Spro-11}). For a subset $A\subseteq\Z^d$, let $\# A$ be the cardinality of $A$, and $\spa A$ be the smallest vector subspace over $\R$ containing $A$.

\subsubsection{Greatest common divisor.} \label{sec:gcd}

For $x,y\in\Z^d$, if $x\neq\bz$, we write $a=\frac{y}{x}$ if there exists $a\in\R$ such that $y=ax$. In particular, $x$ is a \emph{divisor} (\emph{positive divisor}) of $y$, denoted $x|y$, if $\frac{y}{x}\in\Z$ ($\frac{y}{x}\in\N$). For $A\subseteq\Z^d$, $x$ is a {\em common divisor} of $A$, denoted $x|A$, if $x|y$ for all $y\in A$. In particular, $x$ is called a {\em greatest common divisor} (gcd) of $A$, denoted $\gcd(A)$, if $\tilde{x}| x$ for every common divisor $\tilde{x}$ of $A$, and the first non-zero coordinate of $x$ is positive. Hence $\gcd(A)$ is unique, if it exists.
Not all subsets of $\Z^d$ have a gcd (or even a common divisor), e.g., $A=\{(1,2),(2,1)\}$. Indeed, for $A\subseteq\Z^d\setminus\{\bz\}$ non-empty,   $\dim \spa A=1$ if and only if $A$ has a common divisor if and only if $\gcd(A)$ exists (Proposition~\ref{Spro-8}). For a vector $c\in\Z^d$, let $\gcd(c)=\gcd(\{c_j\colon j=1,\ldots,d \})\in\N$.

\subsubsection{Lattice interval.}\label{sec:latint}
For $x,y\in\N^d_0$, denote  the line segment from $x$ to $y$ in $\N^d_0$ by $[x,y]_1$, referred to as the (closed) {\em lattice interval} between $x$ and $y$, and analogously for $[x,y[_1$, etc.
 Moreover, for a non-empty subset $A\subseteq\N^d_0$ on a {\em line} (i.e., $\dim\spa A=1$), if $A$ is not perpendicular to the axis of the first coordinate, then the elements in $A$ are comparable with respect to the  order induced by the first coordinate, denoted $\le_1$.
We use $\min_1 A$ and $\max_1 A$ for the \emph{unique} minimum and maximum of $A$, respectively. If $A$ is countable infinite, then $\max_1 A=\boldsymbol\infty$.
In particular, for $d=1$ and $x,y\in\Z$, $[x,y]_1$ ($[x,y[_1$, respectively, etc.) reduces to the set of  integers from $x$ to $y$.

\subsection{Markov chains}

We first review some basic theory of Markov chains  \cite{N98}.

Let  $Q=(q_{x,y})_{x, y\in\N_0^d}$ be a $Q$-matrix on $\N^d_0$ (not necessarily conservative) \cite{R57}. Let $\cT=\{y-x\colon q_{x,y}>0\ \text{for some}\ x, y\in\N^d_0\}$ be the (possibly infinite) set of jump vectors, and define the \emph{transition rate functions} by
\[\lambda_{\om}\colon \N^d_0\to[0,+\infty), \quad \lambda_\om(x)= q_{x,x+\om},\quad  x\in\N_0^d, \quad \om\in\cT.\]

We say $\omega\in\cT$ is {\em active} in a state $x\in\N^d_0$ if $\lambda_{\om}(x)>0$. In the context of SRNs, $\cT$ is the set of \emph{reaction vectors}. A state $y\in\N_0^d$ is {\em one-step reachable} from $x\in\N^d_0$, denoted $x\rightharpoonup_{\om}y$, if $\om=y-x\in\cT$ and $\om$ is active on $x$. An ordered set of states $\{x^{(j)}\}_{j=1}^m$ for $m>1$ is {\em a path} from $x^{(1)}$ to $x^{(m)}$ if \eqb\label{Eq-9}x^{(1)}\rightharpoonup_{\omega^{(1)}}\cdots\rightharpoonup_{\omega^{(m-1)}}x^{(m)}.\eqe
In particular, if $x^{(1)}=x^{(m)}$, then \eqref{Eq-9} is called a {\em cycle} connecting the states $x^{(i)}$, $i=1,\ldots,m$.
We say a state $y\in\N_0^d$ is {\em reachable} from $x\in\N^d_0$ (or equivalently, $x$ \emph{leads to} $y$) and write $x\rightharpoonup y$
if there exists a path from $x$ to $y$. Hence, in the context of SRNs, $y$ is reachable from $x$ if there is a sequence of reactions converting the molecules in state $x$ into those of $y$. We say $x$ communicates with $y$ and denoted $x \xldownrupharpoon[]{} y$ if both $x\rightharpoonup y$ and $y\rightharpoonup x$ hold. Hence $\xldownrupharpoon[]{}$ defines an equivalence relation on $\N^d_0$, and partitions $\N^d_0$ into {\em communicating classes} (or classes for short). A non-empty subset $E\subseteq\N^d_0$ is {\em closed} if $x\in E\ \text{and}\ x\rightharpoonup y\ \text{implies}\ y\in E$. A set is {\em open} if it is not closed. A state $x$ is {\em absorbing} ({\em escaping}) if $\{x\}$ is a closed (open) class. An absorbing state $x$ is {\em neutral} if $x$ is not reachable from any other state $y\in\N^d_0\setminus\{x\}$, and otherwise, it is {\em trapping}. A non-singleton closed class is a {\em positive irreducible component} (PIC), while a non-singleton open class is a {\em quasi-irreducible component} (QIC). Any singleton class is either an absorbing state (neutral or trapping), or an escaping state, whereas any non-singleton class is either a PIC or a QIC.
Let $\sN$, $\sT$, $\sE$, $\sP$, and $\sQ$ be the (possibly empty) set of all neutral states, trapping states, escaping states, positive irreducible components and quasi-irreducible components for $Q$, respectively. It is easy to verified that every $x\in\N^d_0$ belongs to precisely one of the these sets, and thus all these sets together provide a decomposition of $\N^d_0$.  These sets have intuitive meanings in the context of SRNs, e.g., $\sN\cup\sT$ consists of all states in which no reactions can fire, and $\sE\cup\sP\cup\sQ$ consists of all states in which at least one reaction can fire.

By the definition of these sets, the structure of the ambient space does not depend on the specific form of the transition rate functions but on their \emph{supports} only. The purpose of this paper is to characterize the sets defined above in terms of $\cT$ and the supports of the transition rate functions.

\section{Structure of the state space}\label{sec3}

Let $\cT_{\pm}=\{\omega\in\cT\colon {\rm sgn}(\om_1)=\pm1\}$ be the sets of forward and backward jump vectors, respectively.  Note that $\cT_{\pm}$ is defined in terms of the first coordinate of $\om$. We suppose the following hold.

\medskip
\noindent($\rm\mathbf{A1}$) For every $\om\in\cT$, for $x\in\N^d_0$, $\lambda_{\omega}(x)>0$ implies that $\lambda_{\om}(y)>0$ for all $y\ge x$, $y\in\N^d_0$.
\medskip

\noindent($\rm\mathbf{A2}$) $\cT_+\neq\varnothing$, $\cT_-\neq\varnothing$.
\medskip

We do not require that the CTMCs defined by the $Q$-matrix are unique (i.e., the CTMCs associated with the $Q$-matrix can be \emph{explosive}, or \emph{not minimal}).  Assumption ($\rm\mathbf{A1}$) is a regularity condition on the set of transition functions to evade a CTMC generated by $Q$ to {\em jump sporadically}. It naturally holds for SNRs as well as many population processes modelled as SNRs, where the assumption means a reaction might fire as long as there are enough molecules of the species constituting the reactant \cite{AK11}. If either $\cT_+=\varnothing$ or $\cT_-=\varnothing$, the structural classification is  simpler than under ($\rm\mathbf{A2}$). Indeed, one can derive parallel results from the corresponding results under ($\rm\mathbf{A2}$). Furthermore, by rearrangement of coordinates, ($\rm\mathbf{A2}$) is just required to be fulfilled for one choice of coordinate.

To state the results on  structural classification, we further introduce some notation.
For  $\omega\in\cT$, let $\cI_{\om}$ be the {\em minimal set} of $\supp\lambda_{\om}$.
Let $\cO_{\omega}=\cI_{\omega}+\omega$, $\cI={\cup}_{\omega\in\cT}\cI_{\omega}$ and $\cO={\cup}_{\om\in\cT}\cO_{\omega}$. 
In the context of SRNs,  $\cI_{\om}$ is the set of states in which a reaction with reaction vector $\om$ can fire. $\cI$ refers to the set of \emph{inputs} (or \emph{reactants}), and $\cO$ refers to that of \emph{outputs} (or \emph{products}) \cite{WX20}.

 For ease of notation, we write $\overline{\cI}_{\om}$ for $\overline{\cI_\om}$, etc.
By Proposition \ref{Spro-11}, the set $\cI_{\omega}$ is  non-empty and finite, and by ($\rm\mathbf{A1}$),
$$\overline{\cI}_\omega=\supp\lambda_{\om}, \qquad \overline{\cO}_\omega=\supp\lambda_{\om}+\omega=\overline{\cI}_\omega+\omega,$$
such that the support of a transition rate function is  determined by the corresponding  minimal set. Moreover,
$$\overline{\cI}=\underset{\om\in\cT}{\cup}\overline{\cI}_\omega,\qquad \overline{\cO}=\underset{\om\in\cT}{\cup}\overline{\cO}_\omega.$$

Let $\cF=\{\cI_{\om}\colon \om\in\cT\}$.
In the following, we always take $(\cT,\cF)$ to be associated a $Q$-matrix.

\begin{figure}[h]
\centering
\includegraphics[scale=1]{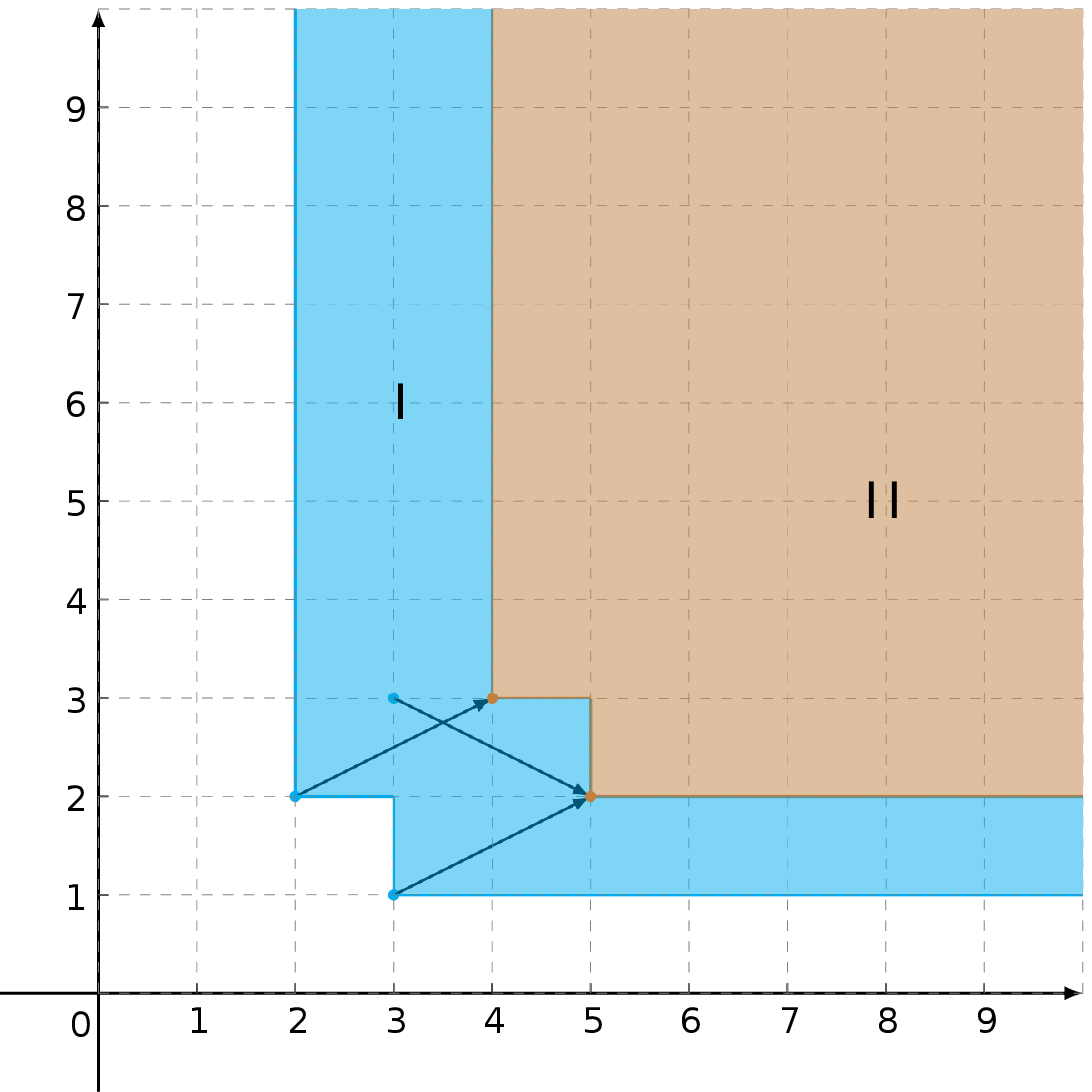}
\caption{Illustration of Example~\ref{ex-1-1}. $\overline{\cI}={\rm I+II}$, $\overline{\cO}={\rm II}$. The elements of the minimal sets $\cI_\om$ are marked with blue dots, those of $\cO_\om$ are brown.
}\label{f5}\end{figure}

\eb\label{ex-1-1} Let  $\Omega=\{(2,1),(2,-1)\}$ and  $\cI_{(2,1)}=\{(2,2),(3,1)\}$, $\cI_{(2,-1)}=\{(3,3)\}$. Then $\cO_{(2,1)}=\{(4,3),(5,2)\}$, $\cO_{(2,-1)}=\{(5,2)\}$. See Figure~\ref{f5} for an illustration.
\ee

\subsection{Characterization of different sets of states}

For  $\om\in\cT$, let
$$\overline{\cI}_{\om}^o=\{x\in\overline{\cI}_{\om}\colon x+\om\rightharpoonup x\}\subseteq \overline{\cI}_{\om},\qquad \cT^o=\{\om\in\cT\colon \overline{\cI}_{\om}^o=\overline{\cI}_{\om}\}\subseteq\cT.$$

The set $\overline{\cI}_{\om}^o$ is the set of states that are reachable from themselves through a cycle of jumps of length at least two, and such that the first jump is $\omega$. The union of $\overline{\cI}_{\om}^o$ over $\omega\in\Omega$ is then the set of states, reachable from themselves through a cycle of jumps of length at least two. It is referred to as \emph{the set of self-reachable states}.

\thmb\label{th-1} Given $(\cT,\cF)$,  assume $\rm(\mathbf{A1})$. Then
$\sN=\N^d_0\setminus(\overline{\cO}\cup \overline{\cI}),$\,  $\sT=\overline{\cO}\setminus \overline{\cI},$\, $\sE=\overline{\cI}\setminus{\bigcup}_{\om\in\cT
}\overline{\cI}_{\om}^o,$\, $\sP\cup\sQ={\bigcup}_{\om\in\cT}\overline{\cI}_{\om}^o,$ and ${\bigcup}_{\om\in\cT}\overline{\cI}_{\om}\setminus\overline{\cI}_{\om}^o\subseteq\sE\cup\sQ$. 
\enb
\item[{\rm(i)}]
$\sN=\N^d_0\setminus\overline{\cI},$ \,$\sP=\overline{\cI},$ \,$\sT=\sE=\sQ=\varnothing$ if and only if $\cT^o=\cT$.
\item[{\rm(ii)}]
$\sP\cup\sQ=\varnothing$ and all states form singleton classes (either neutral, trapping, or escaping), if and only if
$\sum_{\om\in\cT}c_{\om}\om\neq0$ whenever $c_{\om}\in\N_0$, $\om\in\cT$, and $\sum_{\om}c_{\om}\neq0$.
\item[{\rm(iii)}]
$\overline{\cI}\setminus\overline{\cO}\subseteq\sE,$ and  $\sP\cup\sQ\subseteq\overline{\cI}\cap\overline{\cO}$.
\item[{\rm(iv)}]
If for every $x\in{\bigcup}_{\om\in\cT}\overline{\cI}_{\om}^o$, there exists $y\in\N^d_0\setminus{\bigcup}_{\om\in\cT} \overline{\cI}_{\om}^o$ such that $x\rightharpoonup y$, then $\sP=\varnothing$ and $\sQ={\bigcup}_{\om\in\cT}\overline{\cI}_{\om}^o$.
\item[{\rm(v)}]
If there exists $x\in{\bigcup}_{\om\in\cT}\overline{\cI}_{\om}^o$ such that $$x\xldownrupharpoon[]{} y,\quad \text{for all}\quad y\in{\bigcup}_{\om\in\cT}\overline{\cI}_{\om}^o,\ y\neq x,$$
then $\sP={\bigcup}_{\om\in\cT}\overline{\cI}_{\om}^o$, $\sQ=\varnothing$, or $\sQ={\bigcup}_{\om\in\cT}\overline{\cI}_{\om}^o$, $\sP=\varnothing$, and ${\bigcup}_{\om\in\cT}\overline{\cI}_{\om}^o$ consists of a unique non-singleton communicating class.
\ene
Moreover, if there exists $\nu\in\N^d$ such that $\nu\cdot\om=0$ for all $\om\in\cT$, then all communicating classes are finite.
\thme

The proof is given in Section~5. The representations of $\sE$, $\sP\cup\sQ$ as well as $\sE\cup\sQ$ are a bit involved to prove, whereas the rest of the conclusions are fairly straightforward.

We give some intuition  for Theorem~\ref{th-1}. Theorem~\ref{th-1}(i) provides a necessary and sufficient condition for the CTMCs to be  essential, in the sense that all states are in closed classes.  The positive linear independence condition in Theorem~\ref{th-1}(ii)
 $$\sum_{\om\in\cT}c_{\om}\om\neq0\ \textnormal{whenever}\ c_{\om}\in\N_0,\ \om\in\cT,\ \textnormal{and}\ \sum_{\om}c_{\om}\neq0,$$
 ensures that all states form singleton classes, and the CTMCs  have trivial dynamics. Theorem~\ref{th-1}(iii) says that states that can reach other states, but cannot be reached themselves from any other state,  are escaping and thus transient. Theorem~\ref{th-1}(iv)-(v) are concerned with the set of self-reachable states. In (iv), if the set of self-reachable states can escape from any self-reachable state,
 then all non-singleton classes are open. In (v), if the set of self-reachable states is non-empty and 
contained within a single class, 
then there cannot both be PICs and QICs, and the self-reachable states form themselves a class.

In addition, we remark that the conclusion on finiteness of communicating classes in the context of SRNs is known in the literature, c.f. \cite{HJ72,F19,WX20}.

We say  two $Q$-matrices, $Q$ and $\widetilde{Q}$, are {\em structurally equivalent} if for $x,y\in\N^d_0$, $x\rightharpoonup y$ for $Q$ if and only if
 $x\rightharpoonup y$ for $\widetilde Q$.
Hence if $Q$ and $\widetilde{Q}$ are structurally equivalent, then the respective unions of communicating classes for $Q$ and $\widetilde{Q}$ coincide. But the converse is not true, as illustrated below.

Let $(\cT,\cF)$ and $(\widetilde{\cT},\widetilde{\cF})$  be associated with  $Q$ and $\widetilde{Q}$, respectively. We annotate  quantities referring to   $(\widetilde{\cT},\widetilde{\cF})$ with a tilde, for example $\widetilde{\cI}_{\om}$.

\eb\label{ex-1}
(i) Consider the following two reaction networks:
\eqb\label{Eq-1}
\text{S}\ce{<=>[1][2]} 2\text{S}\ce{<=>[7][4]} 3\text{S}\ce{<=>[6][1]} 4\text{S}\ce{->[1]} 5\text{S},\qquad \text{S}\ce{<=>[1][2]}  2\text{S}\ce{<=>[3][1]}  3\text{S}\ce{->[1]}  4\text{S}.
\eqe
The two underlying $Q$-matrices are both on $\N_0$ and associated with $\cT=\{-1,+1\}$ and $\cF$ given by $\cI_1=\{1\}$, $\cI_{-1}=\{2\}$. Hence they are structurally equivalent.
Nonetheless, the associated CTMCs have quite different dynamics: The first network is {\em explosive} while the second is {\em positive recurrent} \cite{ACKK18,XHW20a}.

(ii) Consider $Q$ and $\widetilde{Q}$ associated with $\cT=\{-1,2\}$, $\widetilde{\cT}=\{-1,1\}$, and $\cF=\widetilde{\cF}$ given by $\cI_{-1}=\widetilde{\cI}_{-1}=\{4\}$ and $\cI_2=\widetilde{\cI}_1=\{0\}$. It is easy to verify that $\sE=\{0,1,2\}$, and $\sP=\N\setminus\{1,2\}$ which consists of a unique PIC  for both $Q$ and $\widetilde{Q}$. Nevertheless, $Q$ and $\widetilde{Q}$ are not structurally equivalent, since $0\rightharpoonup1$ for $\widetilde{Q}$ while $0\not\rightharpoonup1$ for $Q$ (however, $0\rightharpoonup2$ is for $\widetilde{Q}$).
\ee

The next result is a direct but non-trivial consequence of Theorem~\ref{th-1}, and provides a criterion for showing  structural equivalence.

\thmb\label{th-0}
Let $Q$ and $\widetilde{Q}$ be associated with $(\cT,\cF)$ and $(\widetilde{\cT},\widetilde{\cF})$, respectively. Assume $\widetilde{\cT}\subseteq\cT$ and $\rm(\mathbf{A1})$ for $\cT$. If
\enb
\item[{\rm(i)}] $\widetilde{\cI}_{\om}=\cI_{\om}$ for all $\om\in\widetilde{\cT}$,
  \item[{\rm(ii)}] for  every $\om\in\cT\setminus\widetilde{\cT}$,  there is a  positive divisor $\widetilde{\om}\in\widetilde\Omega$ of $\om$ such that
 $$x\in\overline{\cI}_{\om}\quad \text{implies}\quad  x+j\widetilde{\om}\in\overline{\cI}_{\widetilde{\om}},\quad \text{for}\quad j=0,\ldots,\tfrac{\om}{\widetilde{\om}}-1.$$
\ene
Then $Q$ and $\widetilde{Q}$ are structurally equivalent.
\thme
\prb

It suffices to show the equivalence of accessibility of one state to another state for $Q$ and $\widetilde{Q}$. For any $x, y\in\N^d_0$, assume $x\rightharpoonup y$ for $Q$. Then there exists $\{\om^{(j)}\}_{j=1}^{m}\subseteq\cT$ such that $x+\sum_{j=1}^{i-1}\om^{(j)}\in\overline{\cI}_{\om^{(i)}}$ for all $i=1,\ldots,m$ and $y=x+\sum_{j=1}^m\om^{(j)}$. Construct $\{\widetilde{\om}^{(j)}\}_{j=1}^{\widetilde{m}}\subseteq\widetilde{\cT}$ from $\{\om^{(j)}\}_{j=1}^{m}$ in the following way:
For every $j=1,\ldots,m$, by conditions (i) and (ii), there exists $k_j\in\N$ and $\widetilde{\om}^{(j)}\in\widetilde{\cT}$ such that $\om^{(j)}=k_j\widetilde{\om}^{(j)}$, $x+\sum_{j=1}^{i-1}\om^{(j)}+\sum_{l=1}^r\widetilde{\om}^{(i)}\in\overline{\cI}_{\widetilde{\om}^{(i)}}$, $\forall i=1,\ldots,m$, $r=0,\ldots,k_i-1$, and $y=x+\sum_{j=1}^{m}\om^{(j)}$. Hence $x\rightharpoonup y$ for $\widetilde{Q}$. The reverse implication is trivial as $\widetilde{\cT}\subseteq\cT$ and condition (i).

By {\rm(i)}-{\rm(iii)}, $x+\sum_{j=1}^{i-1}\om^{(j)}+\sum_{l=1}^{r}\widetilde{\om}^{(j)}\in\overline{\widetilde{\cI}}_{\widetilde{\om}^{(i)}}$,  for $ i=1,\ldots,m$, $r=0,\ldots,k_i-1$, and $y=x+\sum_{j=1}^{m}\om^{(j)}$. Hence $x\rightharpoonup y$ for $\widetilde{Q}$.  The reverse implication is trivial.
\pre

From this result, one can under mild conditions  classify the  ambient space by first reducing the set of jump vectors to a smaller set and then consider the $Q$-matrix associated with this smaller set. Nevertheless, even if $\cT\cap\widetilde{\cT}=\varnothing$, it is still possible for the two $Q$-matrices to be structurally equivalent, as illustrated below.

\eb
Let $Q$ be associated with $\cT=\{-2,1\}$ and $\cF$ given by $\cI_{-2}=\{2\}$, $\cI_1=\{0\}$. Let  $\widetilde{Q}$  be associated with  $\widetilde{\cT}=\{-4,3\}$ and $\widetilde{\cF}$ given by $\widetilde{\cI}_{-4}=\{4\}$, $\widetilde{\cI}_3=\{0\}$.  Although $\cT\cap\widetilde{\cT}=\varnothing$, $\N_0$ is a PIC for both $Q$ and $\widetilde{Q}$.
\ee

Neither of the conditions (i)-(iii) in Theorem~\ref{th-0} can be dropped.

\eb
(i) Let $Q$ be associated with $\cT=\{-1,1\}$ and $\cF$ given by $\cI_{-1}=\{1\}$, $\cI_1=\{0\}$. Let  $\widetilde{Q}$  be associated with $\widetilde{\cT}=\{-1,1\}$ and $\widetilde{\cF}$ given by $\widetilde{\cI}_{-1}=\{1\}$,  $\widetilde{\cI}_1=\{1\}$. Only condition (i) in Theorem~\ref{th-0} fails. Here $\sP=\N_0$ is the unique PIC  for $Q$, while $\sT=\{0\}$  and $\sQ=\N$ for $\widetilde{Q}$. Hence $Q$ and $\widetilde{Q}$ are not structurally equivalent.

(ii) Let $Q$ be associated with  $\cT=\{-2,-1,1\}$ and $\cF$ given by $\cI_{-1}=\cI_{-2}=\{2\}$, $\cI_{1}=\{0\}$.  Let  $\widetilde{Q}$  be associated with
$\widetilde{\cT}=\{-1,1\}$ and $\widetilde{\cF}$ given by $\widetilde{\cI}_{-1}=\{2\}$, $\widetilde{\cI}_{1}=\{0\}$. Only condition (ii) in Theorem~\ref{th-0} fails, since $2\in\overline{\cI}_{-2}$ but $2-1=1\notin\overline{\cI}_{-1}$, and  $\sP=\N_0$  for $Q$, while $\sE=\{0\}$ and $\sP=\N$ for $\widetilde{Q}$. Hence $Q$ and $\widetilde{Q}$  are not structurally equivalent.
\ee

In the setting of population processes, $\sT$ is usually referred to as the {\em extinction set} and $\sT=\varnothing$ implies {\em persistence} of species.
Next, we derive necessary and sufficient conditions in terms of $\cI$ and $\cO$ for $\sT=\varnothing$ and  $\#\sT<\infty$, respectively.

\thmb\label{th-12}
Given $(\cT,\cF)$, assume $\rm(\mathbf{A1})$. Then $\sT=\varnothing$ if and only if  for all  $x\in\cO$, there exists $y\in\cI$ such that $x\ge y$.
\thme

\prb
From Theorem \ref{th-1}, $\sT=\varnothing$ is equivalent to $\overline{\cO}\subseteq\overline{\cI}$. Hence, if $\sT=\varnothing$, choose any $x\in\cO\subseteq\overline{\cO}\subseteq\overline{\cI}$, then there also exists $y\in\cI$ such that $x\ge y$. To see the converse, by definition of $\overline{\cO}$, we find the condition implies $\overline{\cO}\subseteq\overline{\cI}$.
\pre

\thmb\label{th-3}
Given $(\cT,\cF)$, assume $\rm(\mathbf{A1})$.
The following three statements are equivalent:
\begin{itemize}
\item[{\rm (i)}] $\sT$ is finite,
\item[{\rm(ii)}]
$\overline{\widehat{\cO}_j}\subseteq \overline{\widehat{\cI}_j}$, for $j=1,\ldots,d$,
\item[{\rm(iii)}]
for  $j=1,\ldots,d$ and $x\in\cO$, there exists $y\in\cI$ such that  $\widehat x_j\ge \widehat y_j$.
\end{itemize}
\thme
\prb

From Theorem~\ref{th-1}(i),  $\#\sT<\infty$ if and only if $\overline{\cO}\setminus\overline{\cI}$ is bounded, if and only if $(*)$
there exists $M>0$ such that if $x\in\overline{\cO}$ with $x_j>M$ for some $j\in\{1,\ldots,d\}$, then $x\in \overline{\cI}$.

It then suffices to show (iii) is  equivalent to both $(*)$ and (ii).
We first show that $(*)$ is equivalent to (iii), and then (ii) is equivalent to (iii).

 $(*)\Rightarrow$\,(iii). Choose $M$ satisfying $(*)$. Fix $j\in\{1,\ldots,d\}$ and $x\in\cO$. Let $\widetilde{M}=M+1+x_j$, and define $\wx$ by
\[\wx_k=\cab \widetilde{M},\ \text{if}\ k=j,\\ x_k,\ \text{if}\ k\in\{1,\ldots,d\}\setminus\{j\}.\cae\]
Obviously $\wx\ge x$, and $\wx\in \overline{\cO}$ with $\wx_j>M$. Hence, by $(*)$, $\wx\in \overline{\cI}$, i.e., there exists $y\in\cI$ such that $\wx\ge y$, which yields that $x_k=\wx_k\ge y_k$ for all $k\in\{1,\ldots,d\}\setminus\{j\}$, that is, (iii) holds.

 (iii)\,$\Rightarrow(*)$. Set $M=\max_{z\in\cI,1\le i\le d}|z_i|$. Let $x\in \overline{\cO}$ with $x_j>M$ for some  $j\in\{1,\ldots,d\}$. There exists $\ux\in \cO$ such that $x\ge \ux$. By (iii), there exists $y\in\cI$ such that  $x_k\ge \ux_k\ge y_k$ for all $k\in\{1,\ldots,d\}\setminus\{j\}$.  Since $x_j>M\ge y_j$, $x\in \overline{\cI}$, which yields $(*)$.

(ii)\,$\Rightarrow$\,(iii). Fix $j=1,\ldots,d$ and $x\in \cO$. Since
$$\widehat{x}_j\in \widehat{\cO}_j\subseteq \overline{\widehat{\cO}_j}\subseteq\overline{\widehat{\cI}_j},$$
 there exists $y\in \cI$ such that $x_k\ge y_k,$ for all $k\in\{1,\ldots,d\}\setminus\{j\}$, which implies (iii).

(iii)\,$\Rightarrow$\,(ii).
Fix $j=1,\ldots,d$ and $x=(x_1,\ldots,x_{d-1})\in \overline{\widehat{\cO}_j}$. Then there exists $\tilde x\in \cO$ with
$$x_k\ge\left\{\begin{array}{cl} \tilde x_k & \text{if}\quad k=1,\ldots,j-1,\\ \tilde x_{k+1} & \text{if}\quad k=j,\ldots,d-1.\end{array}\right.$$
 By (iii), there exists $y\in \cI$ such that
 $$x_k\ge\left\{\begin{array}{cl}  \tilde x_k\ge y_k &\text{if}\quad k=1,\ldots,j-1,\\ \tilde x_{k+1}\ge y_{k+1} & \text{if}\quad k=j,\ldots,d-1,\end{array}\right.$$
 that is, $x\in\overline{\widehat{\cI}_j}.$
  Since $j=1,\ldots,d$ and $x\in\overline{\widehat{\cO}_j}$ were arbitrary, (ii) holds.
\pre

A necessary condition for  $\#\sT<\infty$ is given below, which is  also sufficient for $d=2$.
\thmb\label{th-4}
Given $(\cT,\cF)$, assume $\rm(\mathbf{A1})$.  If $\sT$ is finite, then \eqb\label{Eq-6}
\min\cO_j\ge\min \cI_j,\quad \text{for}\quad j=1,\ldots,d.
\eqe
\thme

\prb
The necessity follows readily from Theorem \ref{th-3}(iii).
\pre

\thmb\label{th-2}
Given $(\cT,\cF)$, assume $\rm(\mathbf{A1})$. For $d=2$,  $\sT$ is finite if and only if \eqref{Eq-6} holds.
\thme
\prb
By Theorem~\ref{th-4}, it suffices to show that \eqref{Eq-6} is also sufficient for $\#\sT<\infty$.

Let $y^{(j)}\in\cI$ such that $y^{(j)}_j=\min\cI_j$ for $j=1,2$.  From \eqref{Eq-6}, for all $x\in\cO$ and $j=1,2$, we have
\[x_{3-j}\ge\min\cO_{3-j}\ge\min\cI_{3-j}=y^{(3-j)}_{3-j},\]
with $y^{(3-j)}\in\cI$. Since $j\in\{1,2\}$ is arbitrary, Theorem \ref{th-3}(iii) holds.
\pre
Extension of the above result to  dimensions higher than two seems not possible.
\eb
Let $d=3$. Consider $\cT=\{(1,1,1),-(2,2,2)\}$ and $\cF$ given by $\cI_{(1,1,1)}=\{(1,1,3),$ $(1,3,1)\}$, $\cI_{-(2,2,2)}=\{(4,4,4)\}$. Hence $\cI=\{(1,1,3),(1,3,1),$ $(4,4,4)\}$ and $\cO=\{(2,2,4),$ $(2,2,2),$ $(2,4,2)\}$.  For $j=1$ and  $x=(2,2,2)\in\cO$, there exists no $y\in\cI$ such that $x_k\ge y_k$ for $ k=2,3.$  Hence, by Theorem \ref{th-3}, $\#\sT=\infty$.
$$\min\cO_j=2\,>\,\min\cI_j=1,\quad \text{for}\quad j=1,2,3,$$
contradicting  \eqref{Eq-6}.
\ee

\subsection{Classification of states in one dimension}

Let $\sS=\spa\cT$ and note that any CTMC generated by $Q$ with initial state in $\N_0^d$ is confined to an affine  subspace parallel to $\sS$. Define the {\em invariant subspaces} by
\eqb\label{Eq-subspace}\sL_c=\lt(\sS+c\rt)\cap\N_0^d,\quad \text{for}\quad c\in\Z^d,\eqe
which are translationally invariant:
$\sL_c=\sL_{c'}$ whenever $c-c'\in\sS$ and $\sL_c\cap\sL_{c'}=\varnothing$ if $c-c'\not\in\sS$.

In the following we specialize to $\dim\sS=1$. In this case, $\om^*=\gcd(\cT)$ exists (Proposition~\ref{Spro-8})
and  hence $\om^*_1>0$ by definition (Section \ref{sec:gcd}).

Denote the union of the positive and negative minimal sets by $$\cI_{+}=\cup_{\omega\in\cT_{+}}\cI_{\omega}, \quad \cI_{-}=\cup_{\omega\in\cT_{-}}\cI_{\omega}, \quad \cO_{+}=\cup_{\omega\in\cT_{+}}\cO_{\omega},\quad \cO_{-}=\cup_{\omega\in\cT_{-}}\cO_{\omega},$$
and define
 $$\sK_c=\sL_c\cap\Bigl(\bigl(\overline{\cI}_{+}\cap\overline{\cO}_{-}\bigr)\cup\bigl(\overline{\cI}_{-}\cap\overline{\cO}_{+}\bigr)\Bigr),$$
which is independent of the chosen representative due to translational invariance. Let $\uc={\min}_1\sK_c$, and $\bc={\max}_1\sK_c$, such that $\sK_c\subseteq [\uc,\bc]_1$.

We provide some intuition in the context of SRNs. Since $\dim\sS=1$, one might regard $\cT_+$ as the set of reaction vectors corresponding to production events, and similarly, $\cT_-$ as the set of reaction vectors corresponding degradation/consumption events. In the context of generalized BDPs, then $\cT_+$ stores all jumps associated with  birth mechanisms, and $\cT_-$ contains all jumps associated with  death  mechanisms. Here, the set $\overline{\cI}_+$ contains all states in which at least one positive jump is active, and analogously, $\overline{\cO}_-$ contains all states that are one-step reachable from a state in which a negative jump is active. If the CTMC starts in $\sL_c$ then it remains there. Any state in $\sK_c$ can reach another state by a forward jump \emph{and} be reached from some other state in $\sK_c$ by a backward jump, or \emph{vice versa}.

Moreover, define  $\om^{**}=\gcd(\{\om^*_1,\ldots,\om^*_d\})$ and $\sQ_c=\sQ\cap\sL_c$, and   for $k=1,\ldots,\om^{**}$,
$$\Gamma^{(k)}_{c}=\lt(\om^*\Z+\frac{k-1}{\om^{**}}\om^*+c\rt)\cap\sL_c,\quad \sQ^{(k)}_c=\Gamma^{(k)}_c\cap\sQ.$$ Hence, $\sL_c$ is decomposed into $\om^{**}$ different components $\Gamma^{(k)}_c$, and as demonstrated in Theorem~\ref{th-6} below, states of $\sL_c$ can only communicate within the same component. Hence, any PIC or QIC confined to $\sL_c$ coincides with $\sP^{(k)}_c$ or $\sQ^{(k)}_c$, respectively, for some $k$. Other sets like $\sE_c$, $\sE_c^{(k)}$ are defined analogously.
 To distinguish the PICs from the QICs, we define the index sets
$$\Sigma^+_c=\{k\in\{1,\ldots,\om^{**}\}\colon \sQ^{(k)}_c\neq\varnothing\},\qquad \Sigma^-_c= \{1,\ldots,\om^{**}\}\setminus \Sigma^+_c.$$

The following example shows that the above sets might contradict intuition, e.g., a trapping state may be located not on the boundary of $\N^d_0$  but between two escaping states, and $\sK_c$ may  be degenerate or may not coincide with the lattice interval $[c_*,c^*]_1$ (defined in  Section \ref{sec:latint}).

\begin{figure}[h]
\captionsetup[subfigure]{justification=centering, width=6cm, belowskip=-3cm}
\begin{center}
\subfigure[\tiny{$\sP$: Yellow. $\sT$: Red. $\sN$: Black.  $\sE$: Green.}]{
\includegraphics[scale=.5]{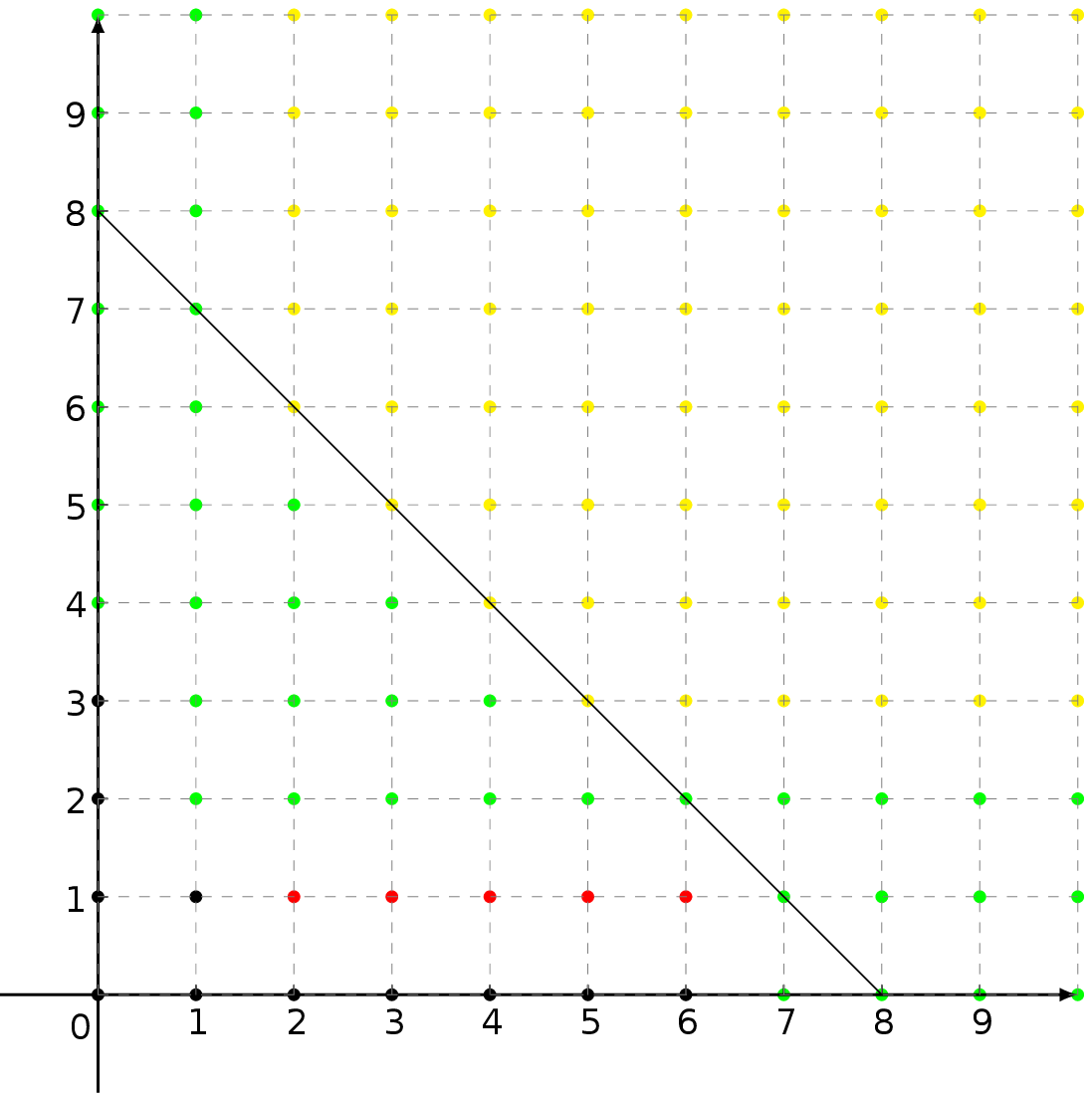}
\hspace{1em}}
\subfigure[{\tiny $\overline{\cI}_+= {\rm I+II+III+V+VI}$, $\overline{\cO}_+={\rm II+III+}$ ${\rm IV+V+VI+VII}$, $\overline{\cI}_-={\rm V+VI+VII +VIII}$, $\overline{\cO}_-$ $={\rm II+V}$, $\big(\overline{\cI}_+\cap\overline{\cO}_-\big)\cup\big(\overline{\cI}_-\cap\overline{\cO}_+\big)={\rm II+V+}$ ${\rm VI+VII}$, $c=(8,0)$.}]
{\includegraphics[scale=.5]{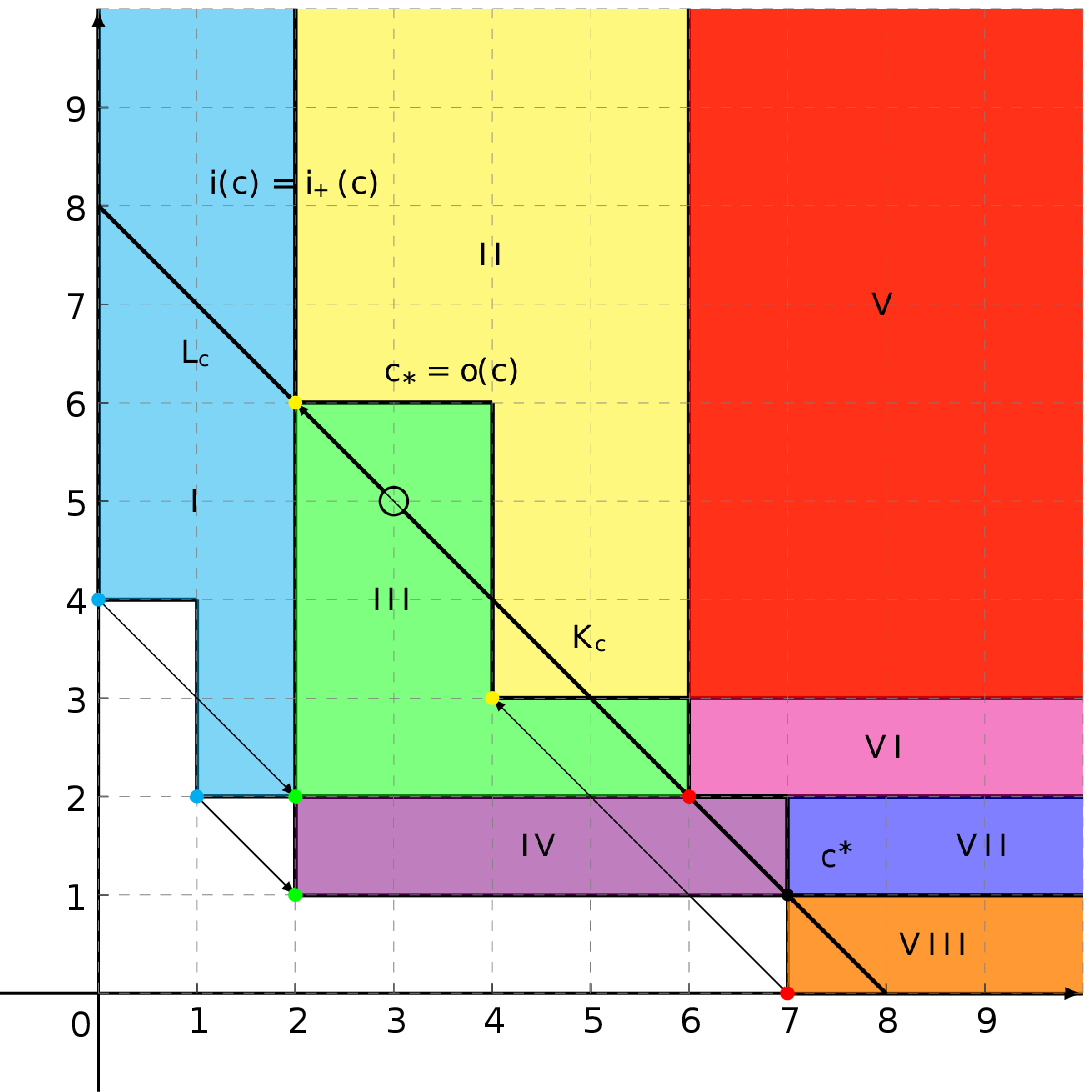}
}
\caption{Illustration of Example \ref{ex-5}.}\label{f2}
\end{center}\end{figure}

\eb\label{ex-5}
Consider $\cT$ and $\cF$ given by $\om^*=(1,-1)$, $\cT_+=\{\om^*,2\om^*\}$, $\cT_-=\{-3\om^*,-4\om^*\}$, $\cI_{\om^*}=\{(1,2)\}$, $\cI_{2\om^*}=\{(0,4)\}$, $\cI_{-3\om^*}=\{(7,0)\}$, and $\cI_{-4\om^*}=\{(6,2)\}$. See Figure~\ref{f2} for an illustration. In this case, $\sT_{(7,0)}=\{(6,1)\}$, $\sE_{(7,0)}=[(0,7),(7,0)]_1\setminus\{(6,1)\}$, $\sK_{(7,0)}=[c_*,c^*]_1=\{(4,3)\}$ and $\sK_{(8,0)}=[(2,6),(7,1)]_1\setminus\{(3,5)\}\subsetneq[c_*,c^*]_1$.
\ee

\thmb\label{th-6} Given $(\cT,\cF)$, assume $\rm(\mathbf{A1})$-$(\rm\mathbf{A2})$ and $\dim\sS=1$.
Then there exists $b\in\N^d_0$ such that for  $c\in\N_0^d+b$, it holds that $\sK_c=[\uc, \bc]_1$, and $\sK_c=\sP_c\cup\sQ_c$ consists of all non-singleton communicating classes  on $\sL_c$, while
$\sL_c\setminus\sK_c$ is the union of singleton communicating classes, composed of
 $$\sN_c= \sL_c\setminus(\overline{\cO}\cup \overline{\cI}),\quad \sT_c=\sL_c\cap \overline{\cO}\setminus \overline{\cI},\quad \sE_c=\overline{\cI}\cap\sL_c\setminus\sK_c.$$
Furthermore,
\begin{itemize}
\item[(i)]  $\sQ^{(k)}_c= \Gamma^{(k)}_c\cap\sK_c$ is a QIC trapped into $\sT^{(k)}_c$ for  $k\in \Sigma^+_c$,
\item[(ii)]  $\sP^{(k)}_c= \Gamma^{(k)}_c\cap\sK_c$ is a PIC for $k\in \Sigma^-_c$.
\end{itemize}
In particular, if $\om^*\in\N_0^d$ then $b=0\in\N_0^d$ suffices.
\thme

The proof is given in Section~6. Theorem~\ref{th-6} provides a classification of states away from the origin. In terms of the input and output sets, $\cI$ and $\cO$, one can identify (within each one-dimensional subspaces $\sL_c$) states of different types. For instance, any state in $\sQ_c^{(k)}$ will eventually be absorbed into the extinction set, and any state in $\sE_c^{(k)}$ will eventually lead to a \emph{stable state}  in the extinction set or in a positive irreducible component.

An escaping state may lead to multiple closed communicating classes, e.g., both a trapping state and a PIC.

\eb\label{ex:closed-escape}
Consider $\cT=\{(1,-1),(-1,1)\}$ and $\cF$ given by $\cI_{(1,-1)}=\{(1,2),(2,1)\}$ and $\cI_{(-1,1)}=\{(1,2),(3,0)\}$. Here $\om^*=(1,-1)$. For  $k\in\N\setminus\{1,2\}$, let $c(k)=(k,0)$. Then $\sL_{c(k)}=\{(j,k-j)\}_{j=0}^k$, $\sT_{c(k)}=\{(0,k)\}$, $\sE_{c(k)}=\{(1,k-1)\}$, and $\sP_{c(k)}=\sL_{c(k)}\setminus\{(0,k),(1,k-1)\}$ which consists of a unique PIC on $\sL_{c(k)}$. Hence $\om^{**}=1$, $\Sigma^+_{c(k)}=\varnothing$ and $\Sigma^-_{c(k)}=\{1\}$, and the unique escaping state leads to both a trapping state and a PIC on $\sL_{c(k)}$.
\ee

To characterize the ambient space further, we define
$$\si(c)={\min}_1\sL_c\cap\overline{\cI},\quad \si_+(c)={\min}_1\sL_c\cap\overline{\cI}_+,\quad \so(c)={\min}_1\sL_c\cap\overline{\cO},\quad \so_-(c)={\min}_1\sL_c\cap\overline{\cO}_-,$$
 for $c\in\N_0^d$. By definition, $\uc\ge_1\si(c),\so(c)$. The result below provides detailed characterization of the relevant sets for $\om^*\in\N_0^d$  and excludes the possibility for an escaping state to lead to both a PIC and a trapping state as shown in Example \ref{ex:closed-escape}.

\cob\label{co-1}
Given $(\cT,\cF)$, assume $\rm(\mathbf{A1})$-$\rm(\mathbf{A2})$, $\dim\sS=1$, and $\om^*\in\N_0^d$. Then for $c\in\N_0^d$, we have $\uc={\max}_1\lt\{\si_+(c),\so_-(c)\rt\}$, $\bc=\boldsymbol{+\infty}$, and
$$\sN_c=\lt[{\min}_1 \sL_c,{\min}_1\{\si(c),\so(c)\}\rt[_1,\quad \sT_c=[\so(c),\si(c)[_1,\quad \sE_c=\lt[\si(c),\uc\rt[_1.$$
These are all finite sets.
Moreover,
 $$\Sigma_c^+=\Bigl\{1+\frac{\om^{**}(v-c)}{\om^*}-\Bigl\lf\frac{v-c}{\om^*}\Bigr\rf\om^{**}\colon v\in[\so(c),\si(c)[_1\Bigr\},$$
 $$\Sigma_c^-=\Bigl\{1+\frac{\om^{**}(v-c)}{\om^*}-\Bigl\lf\frac{v-c}{\om^*}\Bigr\rf\om^{**}\colon v\in[\si(c),\so(c)+\om^*[_1\Bigr\},$$
 and
 \begin{align*}
 \Gamma_c^{(k)}\cap(\om^*\N_0+\uc)&=\om^*\Bigl(\N_0+\Bigl\lc\frac{c_*-c-\frac{k-1}{\om^{**}}\om^{*}}{\om^*}\Bigr\rc\Bigr)+\frac{k-1}{\om^{**}}\om^{*}+c \\
 &=\cab\sP_c^{(k)},\quad\text{if}\ k\in\Sigma_c^-,\\ \sQ_c^{(k)},\quad\text{if}\ k\in\Sigma_c^+.\cae
 \end{align*}
Hence $\# \Sigma_c^+={\min}\lt\{\om^{**},\max\{0,\tfrac{\om^{**}(\si(c)-\so(c))}{\om^*}\}\rt\}$. Furthermore, if  $\Sigma^+_c\neq\varnothing$, then $\sT_c\neq\varnothing$, which moreover implies that any escaping state  never reaches a state in $\sP_c$ or $\sQ_c$, but only another escaping state or a trapping state.
\coe
The proof is given in Section~7. The number of PICs and QICs has been obtained for one-species SRNs \cite[Corollary~5.3]{HW20}. Here, we characterize these components in more detail in a general context beyond SRNs.

The following example illustrates that a QIC may lead to multiple trapping states, and an escaping state may also lead to multiple escaping or trapping  states, even if $\om^*\in\N_0^d$.
\eb Let $d=1$. Consider
$\cT=\{-2,-1,1\}$ and $\cF$ given by $\cI_{-2}=\cI_{-1}=\{2\}$, and $\cI_1=\{5\}$.
The flow chart for the  state space is as follows:
\[
\begin{tikzpicture}[node distance=3em, auto, scale=1]
 \tikzset{
    pil/.style={
           ->,
           shorten <=2pt,
           shorten >=2pt,}
}
\node[] (a0) {0};
\node[right=of a0] (a1) {1};
\node[right=of a1] (a2) {2};
\node[right=of a2] (a3) {3};
\node[right=of a3] (a4) {4};
\node[right=of a4] (a5) {5};
\node[right=of a5] (a6) {6};
\node[right=of a6] (a7) {7};
\node[right=of a7] (a8) {$\cdots$};
\node[left=of a6] (aa5) {} edge[pil, black, bend left=20] (a6);
\node[left=of a7] (aa6) {} edge[pil, black, bend left=20] (a7);
\node[left=of a8] (aa7) {} edge[pil, black, bend left=15] (a8);
\node[right=of a7] (aa8) {} edge[pil, black, bend right=0] (a7);
\node[right=of a6] (aa7) {} edge[pil, black, bend right=0] (a6);
\node[right=of a5] (aa6) {} edge[pil, black, bend right=0] (a5);
\node[right=of a4] (aa5) {} edge[pil, black, bend right=0] (a4);
\node[right=of a3] (aa4) {} edge[pil, black, bend right=0] (a3);
\node[right=of a2] (aa3) {} edge[pil, black, bend right=0] (a2);
\node[right=of a1] (aa2) {} edge[pil, black, bend right=0] (a1);
\node[right=of a7] (aa8) {} edge[pil, black, bend left=30] (a6);
\node[right=of a6] (aa7) {} edge[pil, black, bend left=30] (a5);
\node[right=of a5] (aa6) {} edge[pil, black, bend left=30] (a4);
\node[right=of a4] (aa5) {} edge[pil, black, bend left=30] (a3);
\node[right=of a3] (aa4) {} edge[pil, black, bend left=30] (a2);
\node[right=of a2] (aa3) {} edge[pil, black, bend left=30] (a1);
\node[right=of a1] (aa2) {} edge[pil, black, bend left=30] (a0);
\end{tikzpicture}
\]
For all $c\in\N_0$, $\so(c)=0$ and $\si(c)=2$. Then $\sQ_c=\N\setminus\{1,2,3,4\}$ consists of a unique QIC, $\sE_c=\{2,3,4\}$ and $\sT_c=\{0,1\}$. The QIC eventually leads to both escaping states, and the state escaping $4$ leads to the escaping states $2$ and $3$. Eventually both $\sQ_c$ and $\sE_c$ lead to  the two trapping states $0$ and $1$.
\ee

\begin{figure}[h]
\captionsetup[subfigure]{justification=centering,belowskip=-2cm}
\begin{center}
\subfigure[\tiny Example \ref{ex-3}(i).]{
\includegraphics[scale=.5]{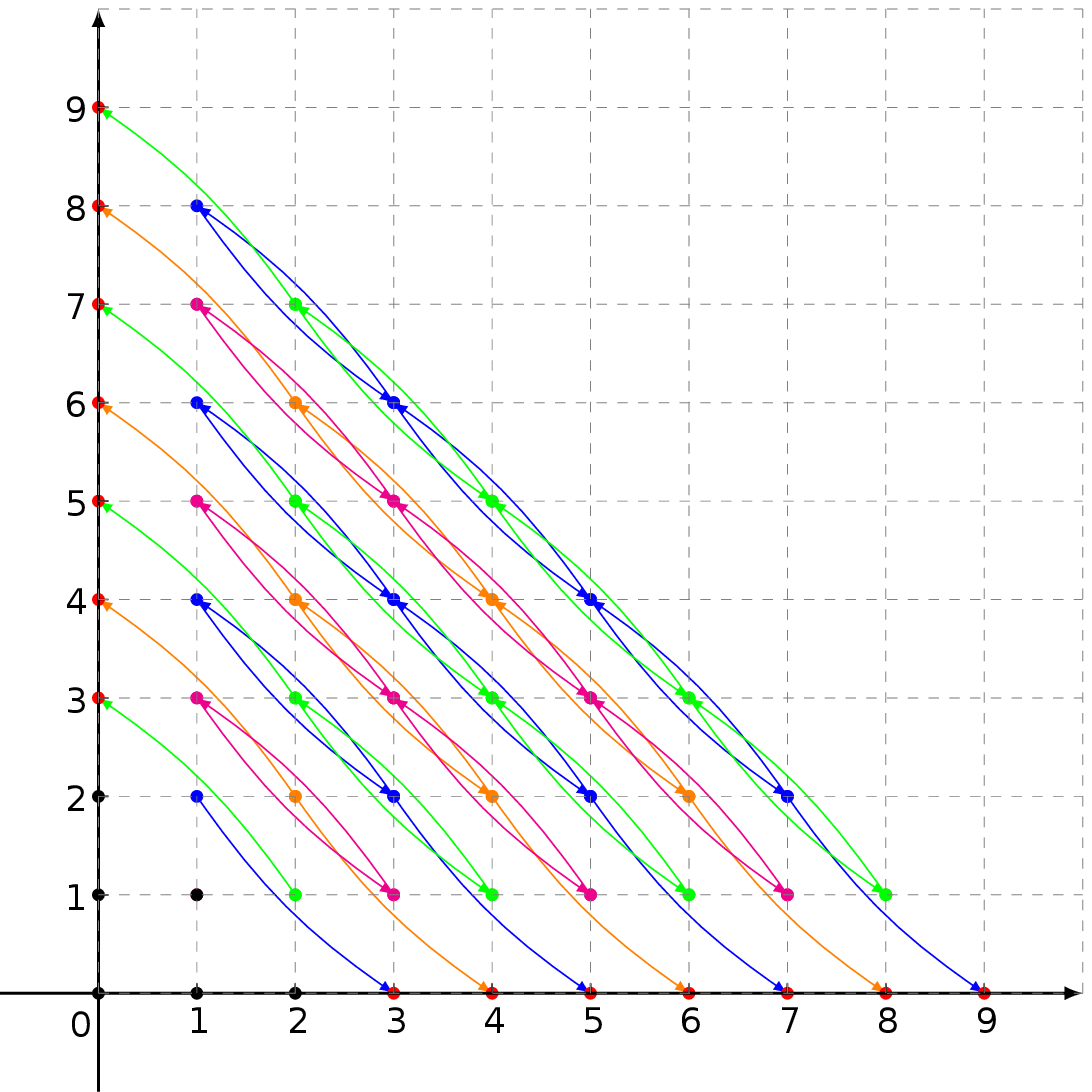}
\hspace{1em}}
\subfigure[\tiny Example \ref{ex-3}(ii).]{
\includegraphics[scale=.5]{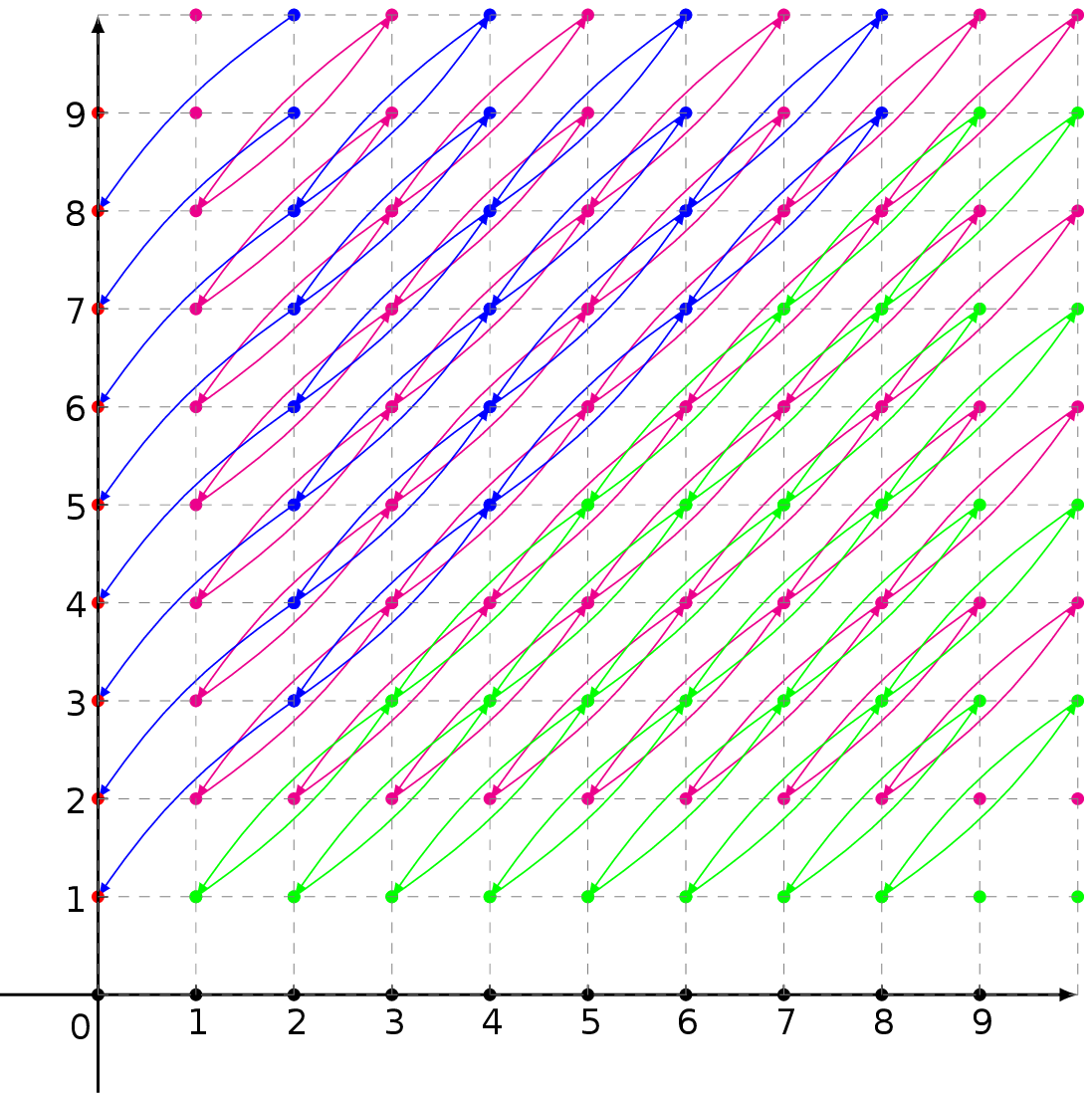}
}
\caption{Illustration of Example \ref{ex-3}. Coexistence of PICs and QICs. (a) $\sN$: Black. $\sT$: Red. QIC (together with the escaping states $(1,2)$ and $(2,1)$, trapped into a single trapping state): Green+blue (different colors for different components). QIC (trapped into both trapping states): Orange. PIC: Magenta. Here PICs and QICs are finite. (b) $\sN$: Black. $\sT$: Red. QIC (trapped into a single  state): Blue.  PIC: Magenta+green. Here PICs and QICs are infinite. }\label{f3}
\end{center}\end{figure}

The following example illustrate how Corollary~\ref{co-1} can be used to classify the ambient space when $\dim\sS=1$.

\eb\label{ex-3} We consider two examples on $\N^2_0$ where PICs and QICs  coexist on the same invariant subspace: {\rm(i)}   $\om^*=(2,-2)$, $\cT=\{\om^*,-\om^*\}$ and $\cF$ given by $\cI_{\om^*}=\{(1,2)\}$, $\cI_{-\om^*}=\{(2,1)\}$, and {\rm(ii)}  $\om^*=(2,2)$, $\cT=\{\om^*,-\om^*\}$ and $\cF$ given by $\cI_{\om^*}=\{(1,1)\}$, $\cI_{-\om^*}=\{(2,3)\}$.  The structure of $\N^2_0$ is illustrated in Figure~\ref{f3}.  We further calculate all  communicating classes for the second $Q$-matrix. Note that $\om^{**}=2$.

(i) Let $c(k)=(0,k)$ for $k\in\N\setminus\{1,2,3,4\}$. Then $\sL_{c(k)}=(\frac{1}{2}\om^*\N_0+c(k))\cap\N^2_0=\{(z_1,z_2)\in\N^2_0\colon z_1+z_2=k\}$. By definition, $\sL_c=[c(k),(k,0)]_1$, $\sK_{c(k)}=[c_*,c^*]_1$ with $c_*=(1,k-1)$ and $c^*=(k-1,1)$,  $\sN_{c(k)}=\sE_{c(k)}=\varnothing$, $\sT_{c(k)}=\{c(k),(k,0)\}$. Moreover, by straightforward calculation, when $k$ is even, $\sQ_{c(k)}=\{(2,k-2),(4,k-4),\ldots,(k-2,2)\}$ consists a single QIC and $\sP_{c(k)}=\{(1,k-1),(3,k-3),\ldots,(k-1,1)\}$ consists a single PIC; when $k$ is odd, $\sQ^1_{c(k)}=\{(2,k-2),(4,k-4),\ldots,(k-2,2)\}$ and $\sQ_{c(k)}^2=\{(1,k-1),(3,k-3),\ldots,(k-1,1)\}$ and $\sP_{c(k)}=\varnothing$. See Figure~\ref{f3}(a).

(ii) Let $c(k)=(k,0)$ for $k\in\N_0$. Then $\sL_{c(k)}=(1,1)\N_0+c(k)$, $\so_-(c(k))=\so(c(k))=\si_+(c(k))=\si(c(k))=(k+1,1)$. By Corollary~\ref{co-1}, $c_*(k)=(k+1,1)$, $\sK_{c(k)}=(1,1)\N+c(k)$, $\sN_{c(k)}=\{c(k)\}$, $\sE_{c(k)}=\sT_{c(k)}=\varnothing$, $\Sigma_{c(k)}^+=\varnothing$, and $\Sigma_{c(k)}^-=\{1,2\}$. Hence $\sP^{(1)}_{c(k)}=(2,2)\N+c(k)$ and $\sP^{(2)}_{c(k)}=(2,2)\N_0+(k+1,1)$. Similar conclusions can be derived for $c(k)=(0,k)$ for $k\in\N$. See Figure~\ref{f3}(b).
\ee

\section{Applications}\label{sec4}

\subsection{Stochastic reaction networks}
SRNs are used to describe interactions of constituent molecular species, though the area of application extends beyond (bio)chemistry \cite{G83,PCMV15}. In this section, we apply the main results from Section~\ref{sec3} to some examples of SRNs to see how diverse the structure of the ambient space $\N_0^d$ is.

An SRN  is a class of CTMCs on $\N_0^d$ given by a reaction graph. A reaction graph is a directed edge-labelled graph where each edge is a reaction $y\ce{->[$\ka_{y\to y'}$]}y'$, $y,\ y'\in\N^d_0$, and the label a positive \emph{reaction rate constant}.
We consider SRNs with \emph{mass-action kinetics} \cite{AK11}, that is, to each reaction is associated a \emph{reaction rate function}
$$x\mapsto\ka_{y\to y'}x^{\underline{y}},\quad x\in\N^d_0,$$ expressing the propensity of a reaction to occur, where $x^{\underline{y}}=\prod_{j=1}^d\prod_{i=0}^{y_j-1}(x_j-i)$ denotes the descending factorial. Hence the transition rate functions for the underlying CTMCs associated with an SRN is the sum of reaction rate functions contributing to the same jump vector:
\[\lambda_{\om}(x)=\sum_{y\to y',\, y'-y=\om}\ka_{y\to y'}x^{\underline{y}},\quad \om\in\cT,\]
where $\cT$ is the set of all reaction vectors $y'-y$ without multiplicity.

An SRN is \emph{essential} if $\N^d_0$ can be decomposed into a disjoint union of closed communicating classes (i.e., $\N^d_0=\sN\cup\sP$). In particular, essential SRNs are \emph{persistent} (in the sense of an empty extinction set).

Structural equivalence of stochastic reaction networks does not necessarily imply \emph{dynamical equivalence} thereof, as illustrated below, and already alluded to in Example~\ref{ex-1}.

\eb
Consider the following two reaction networks:
\[1)\quad \varnothing \ce{<=>[$\ka_1$][$\ka_2$]} \tS,\qquad \text{and}\qquad  2)\quad\varnothing \ce{<=>[$\ka_1$][$\ka_2$]} \tS,\quad 2\tS \ce{->[$\ka_3$]}3\tS,\]
with propensities  given by
$$\lambda^1_{1}(x)=\kappa_1,\quad \lambda^1_{-1}(x)=\kappa_2 x,\qquad \text{and}$$
$$\lambda^2_{1}(x)=\kappa_1+\kappa_3 x(x-1),\quad \lambda^2_{-1}(x)=\kappa_2 x, $$
respectively.
For both SRNs,  $\cT=\widetilde{\cT}=\{-1,1\}$ and $\cF=\widetilde{\cF}$ with  $\cI_{-1}=\widetilde{\cI}_{-1}=\{1\}$ and $\cI_1=\widetilde{\cI}_1=\{0\}$.  By Theorem~\ref{th-0}, the two reaction networks are structurally equivalent. Nevertheless, the first SRN is
 is positive recurrent and admits an exponentially ergodic stationary distribution on $\N_0$. In contrast, the second reaction network is   explosive a.s. for any initial state \cite{XHW20a}.
\ee

\eb
Consider the following SRN:
\[\tS_1+\tS_2\ce{->[$\kappa_1$]}2\tS_2,\quad 2\tS_1+\tS_2\ce{->[$\kappa_2$]}3\tS_1+2\tS_2,\]
$$\lambda_{(-1,1)}(x_1,x_2)=\kappa_1x_1x_2,\quad \lambda_{(1,1)}(x_1,x_2)=\kappa_2x_1(x_1-1)x_2.$$
For this reaction network, $\cT=\{(-1,1),(1,1)\}$ and $\cF$ is given by $\cI_{(-1,1)}=\{(1,1)\}$, $\cI_{(1,1)}=\{(2,1)\}$. Hence $\cI=\{(1,1),(2,1)\}$, $\cO=\{(0,2),(3,2)\}$. Since $(0,2)\not\ge(1,1)$, $(0,2)\not\ge(2,1)$, by Theorem~\ref{th-12}, the extinction set $\sT\neq\varnothing$. Moreover, $\min\cO_1=0<\min\cI_1=1$, by Theorem~\ref{th-2}, $\sT$ is countably  infinite.
\ee

\eb Recall the three SRNs in the Introduction:
\[
\begin{tikzpicture}[node distance=2.5em, auto]
 \tikzset{
    pil/.style={
           ->,
           shorten <=2pt,
           shorten >=2pt,}
}
 \node[] (a) {};
  \node[right=1.2em of a] (n1) {};
  \node[above=-.5em of n1] (m1) {$\ka_1$};
\node[right=of a] (b) {2\tS};
\node[right=1.3em of b] (n2) {};
  \node[above=.9em of b] (m3) {$\ka_3$};
  \node[left=of b] (aa) {\tS} edge[pil, black, bend left=0] (b);
\node[right=of b] (c) {3\tS;};
\node[right=of b] (cc) {} edge[pil, black, bend right=40] (aa);
\node[right=.8em of c] (d) {};
  \node[right=.95em of d] (n11) {};
  \node[above=-.5em of n11] (m11) {$\ka_1$};
\node[right=of d] (e) {2\tS};
\node[right=.8em of e] (n22) {};
  \node[above=-.5em of n22] (m22) {$\ka_2$};
  \node[above=.9em of e] (m33) {$\ka_3$};
  \node[left=of e] (dd) {\tS} edge[pil, black, bend left=0] (e);
\node[right=of e] (f) {3\tS;};
 \node[left=of f] (ee) {} edge[pil, black, bend left=0] (f);
\node[right=of e] (ff) {} edge[pil, black, bend right=40] (dd);
\node[right=.8em of f] (g) {};
  \node[right=.9em of g] (n111) {};
  \node[above=-.5em of n111] (m111) {$\ka_1$};
\node[right=of g] (h) {2\tS};
\node[right=.8em of h] (n222) {};
  \node[above=-.5em of n222] (m222) {$\ka_2$};
  \node[above=.9em of h] (m333) {$\ka_3$};
  \node[left=of h] (gg) {\tS} edge[pil, black, bend left=0] (h);
\node[right=of h] (i) {3\tS};
  \node[right=1em of i] (n444) {};
  \node[above=-.5em of n444] (m444) {$\ka_4$};
 \node[left=of i] (hh) {} edge[pil, black, bend left=0] (i);
\node[right=of h] (ii) {} edge[pil, black, bend right=40] (gg);
\node[right=of i] (j) {4\tS.};
 \node[left=of j] (ii) {} edge[pil, black, bend left=0] (j);
\end{tikzpicture}\]
 By Theorem~\ref{th-0}, the three SRNs are structurally equivalent. The second SRN is \emph{weakly reversible} (whose reaction graph is strongly connected), but the other two SRNs are not. Since weakly reversible SRNs are essential \cite{PCK14}, the above three SRNs are all essential and thus persistent.
\ee

\subsection{ Lotka-Volterra system}
Consider the Lotka-Volterra system represented as a reaction network \cite{HW20}:
\[\tS_1\ce{->[$\ka_1$]}2\tS_1,\quad \tS_1+\tS_2\ce{->[$\ka_2$]}2\tS_2,\quad \tS_2\ce{->[$\ka_3$]}\varnothing,\]
where
$$\lambda_{(1,0)}(x_1,x_2)=\kappa_1 x_1,\quad \lambda_{(0,1)}(x_1,x_2)=\kappa_2 x_1x_2,\quad \lambda_{(0,-1)}(x_1,x_2)=\kappa_3 x_2.$$
Properties of the Lotka-Volterra system have received previous attention \cite{HW20,Klebaner}. Here we formally decompose the ambient space $\N^2_0$ into disjoint classes.

We have, $\cT=\{(1,0),(-1,1),(0,-1)\}$ and $\cF$ is given by $\cI_{(1,0)}=\{(1,0)\}$, $\cI_{(-1,1)}=\{(1,1)\}$, $\cI_{(0,-1)}=\{(0,1)\}$. Moreover, $\cO_{(1,0)}=\{(2,0)\}$, $\cO_{(-1,1)}=\{(0,2)\}$, $\cO_{(0,-1)}=\{(0,0)\}$, $\overline{\cI}=\N_0^2\setminus\{(0,0)\}$,  $\overline{\cO}=\N_0^2$. Hence
 $$\overline{\cI}_{(1,0)}=\N\times\N_0,\quad  \overline{\cI}_{(-1,1)}= \N\times\N, \quad  \overline{\cI}_{(0,-1)}=\N_0\times\N.$$
Furthermore, since $\sum_{\om\in\cT}\om=\bz$, one can  verify that $\overline{\cI}^o_{(1,0)}=\N^2\subsetneq\overline{\cI}_{(1,0)}$, $\overline{\cI}^o_{(-1,1)}=(\N\setminus\{1\})\times\N\subsetneq\overline{\cI}_{(-1,1)}$, $\overline{\cI}^o_{(0,-1)}=\N\times(\N\setminus\{1\})\subsetneq\overline{\cI}_{(0,-1)}$. Hence $\cT^o=\varnothing$. By Theorem~\ref{th-1},
$$\sN=\varnothing,\quad \sT=\overline{\cO}\setminus\overline{\cI}=\{(0,0)\},\quad \sP\cup\sQ=\N^2,\quad
\sE=\overline{\cI}=(\{0\}\times\N)\cup(\N\times\{0\}).$$

    Moreover, for every $(x,y)\in\N^2=\sP\cup\sQ\subseteq\overline{\cI}_{(0,-1)}$, we have $(x,y)\rightharpoonup(x,0)\in\N_0\times\{0\}\subsetneq\sT\cup\sE$, by repetition of the jump vector $(0,-1)$. Hence, by Theorem~\ref{th-1}(iv), $$\sP=\varnothing,\qquad \sQ=\N^2.$$
   Finally, we show  $\sQ=\N^2$ consists of a unique QIC. For any $(m,n)\in\N^2\setminus\{(1,1)\}$, $(1,m+n-1)\rightharpoonup(1,1)$ by repetition of the jump vector $(0,-1)$. In addition, if $m\neq1$, then $(m,n)\rightharpoonup(1,m+n-1)$ by repetition of the jump vector $(-1,1)$.  Hence $(m,n)\rightharpoonup(1,1)$. Similarly, we have $(1,1)\rightharpoonup(m+n-1,1)$ by repetition of the jump vector $(1,0)$, and $(m+n-1,1)\rightharpoonup(m,n)$ by repetition of the jump vector $(-1,1)$, provided $n\neq 1$. In sum, $(m,n)\xldownrupharpoon[]{} (1,1)$. The conclusion now follows from Theorem~\ref{th-1}(v).

\subsection{The EnvZ-OmpR system}

Consider the EnvZ-OmpR signaling pathway in {\em Escherichia coli} \cite{SF10}:
\[\tS_1\ce{<=>[$\ka_1$][$\ka_2$]}\tS_2\ce{<=>[$\ka_3$][$\ka_4$]}\tS_3\ce{->[$\ka_5$]}\tS_4\qquad \tS_4+\tS_5\ce{<=>[$\ka_6$][$\ka_7$]}\tS_6\ce{->[$\ka_8$]}\tS_2+\tS_7\]
\[\tS_3+\tS_7\ce{<=>[$\ka_9$][$\ka_{10}$]}\tS_8\ce{->[$\ka_{11}$]}\tS_3+\tS_5\qquad \tS_1+\tS_7\ce{<=>[$\ka_{12}$][$\ka_{13}$]}\tS_9\ce{->[$\ka_{14}$]}\tS_1+\tS_5\]

We use generic names for the species ($S_1,\ldots,S_9$) for simplicity. These correspond, however, to specific molecular biology compounds in the EnvZ-OmpR bacterial system.

Let $e_i\in\N^9_0$ denote the unit vector with 1 in the $i$-th coordinate and 0 in the remaining coordinates, for $i=1,\cdots,9$. Then, the 14 \emph{ordered} (indicated by the indices of reaction rate constants $\ka_i$) reactions contribute to 14 different jump vectors denoted $\om^{(i)}$, for $i=1,\ldots,14$. For instance, $\om^{(1)}=e_2-e_1$, derived from the first reaction. Hence $\cT=\{\om^{(i)}\}_{i=1}^{14}$.
It is straightforward to verify that $$\overline{\cI}=\{x\in\N^9_0\colon \text{either}\  x_1+x_2+x_3+x_6+x_8+x_9\ge1,\ \text{or}\  x_4\ge1,  x_5\ge1\ \text{holds}\},$$
 which implies that $\sN\cup\sT=\N^d_0\setminus\overline{\cI}=\{ie_4+je_5+ke_7\colon ij=0,\ i,j, k\in\N_0\}$. Moreover, $\overline{\cO}\setminus\overline{\cI}=\{e_4\}$, which implies that $\overline{\cO}\setminus\overline{\cI}=\overline{\{e_4\}}\setminus\overline{\cI}$. Hence $$\sT=\overline{\cO}\setminus\overline{\cI}=
\lt\{ie_4+je_7\colon i\in\N,\ j\in\N_0\rt\}\neq\varnothing,\quad \sN=\lt\{ie_5+je_7\colon i, j\in\N_0\rt\}.$$

Let $\widetilde{\cT}=\cT\setminus\{\om^{(i)}\colon i=5,8,11,14\}$. Then $\cup_{\om\in\widetilde{\cT}}\cI_{\om}=\cup_{\om\in\cT}\cI_{\om}$. Moreover, for $\om\in\widetilde{\cT}$,  we have $-\om\in\cT$, and $\cI_{\om}=\cO_{-\om}$, which implies that $\overline{\cI}_{\om}^o=\overline{\cI}_{\om}$.
Since $e_4\in\sT$, we have $e_4\not\rightharpoonup e_3$, and thus $e_3\in\overline{\cI}_{\om^{(5)}}\setminus\overline{\cI}_{\om^{(5)}}^o$.
However, $\overline{\cI}_{\om^{(8)}}=\overline{\cI}_{\om^{(8)}}^o$, since there is a path
$e_2+e_7\rightharpoonup e_6$ (applying reactions 3, 9, 11, 5, 6 in that order).
Similarly, $\overline{\cI}_{\om^{(11)}}=\overline{\cI}_{\om^{(11)}}^o$ and $\overline{\cI}_{\om^{(14)}}=\overline{\cI}_{\om^{(14)}}^o$, since there is a path
$e_3+e_5\rightharpoonup e_8$ (applying reactions 5, 6, 8, 3, 9) and $e_1+e_5\rightharpoonup e_9$ (applying reactions 1, 3, 5, 6, 8, 2, 12).
This shows $\overline{\cI}_{\om}^o=\overline{\cI}_{\om}$ if and only if $\om\in\cT\setminus\{\om^{(5)}\}$, and hence $\cT^o=\cT\setminus\{\om^{(5)}\}$. Since $\cT^o\subseteq\widetilde{\cT}\subseteq\cT$ and $\cup_{\om\in\widetilde{\cT}}\cI_{\om}=\cup_{\om\in\cT}\cI_{\om}$, we have $\bigcup_{\om\in\cT}\overline{\cI}_{\om}^o=\overline{\cI}$.
By Theorem~\ref{th-1}, we have $\sE=\varnothing,$ and
$$\sP\cup\sQ=\{x\in\N^9_0\colon \text{either}\ x_1+x_2+x_3+x_6+x_8+x_9\ge1,\ \text{or}\  x_4\ge 1, x_5\ge1\ \text{holds}\}.$$

Finally, we show $\sP=\varnothing$, and $\sQ=\overline{\cI}$. By Theorem~\ref{th-1}, it suffices to show that for every $x\in\overline{\cI}$, there exists $y\in\sT\cup\sN=\{ie_4+je_5+ke_7\colon ij=0,\ i,j,k\in\N_0\}$ such that $x\rightharpoonup y$. Similarly, we can use respective reactions to find the following paths:
$$e_i\rightharpoonup e_4,\quad \text{for}\ i=1, 2, 3,\quad e_i\rightharpoonup e_4+e_5,\quad \text{for}\ i=6, 8, 9,\quad e_4+e_5\rightharpoonup e_4+e_7.$$
 Using these paths, any state in $\overline{\cI}$ with non-zero entries in the  coordinates $1, 2, 3, 6, 8, 9$ leads to one state with zero entries in these coordinates, and any state therein with non-zero entries in both coordinates $4, 5$ leads to another state with zero entry in coordinate $5$. In this way, we show that any state in $\overline{\cI}$ leads to a state in $\{ie_4+ke_7\colon i,j,k\in\N_0\}\subseteq\sT\cup\sN$.
Hence we conclude that $\N^9_0$ is decomposed into non-empty QICs and single closed classes. Since, $\nu=(1, 1, 1, 1, 1, 2, 1, 2, 2)$ is orthogonal to all jump vectors in $\Omega$, then all QICs are finite.

This further indicates that
any CTMC associated with this SRN has \emph{certain absorption} and admits a unique ergodic QSD supported on a QIC \cite{CMM13}. In \cite{AEJ14}, an analysis of the ambient space is conducted and it is concluded that
\[(\sP\cup\sT\cup\sN)\cap(\cup_{\om\in\widetilde{\cT}}\overline{\cI}_{\om})=\varnothing.\]
(in the terminology of this paper).
However, this fails to conclude $\sP=\varnothing$ directly, and also  $\sE=\varnothing$. Thus, our approach, although tedious in the specific example,  provide extra information of the system beyond the literature.

\subsection{Extended class of branching processes}

Consider an extended class of branching processes with $Q$-matrix $Q=(q_{x,y})_{x,y\in\N_0}$:
\eqb\label{Eq-5}q_{x,y}=\left\{\begin{array}{cl}
r(x)\mu(y-x+1), &\quad  \text{if}\quad  y\ge x-1\ge0\quad \text{and}\quad y\neq x,\\
 -r(x)(1-\mu(1)), &\quad  \text{if}\quad y=x\ge1,\\
 q_{0,y}, & \quad \text{if}\quad y>x=0,\\
 -q_0, &\quad \text{if}\quad y=x=0,\\
 0, &\quad  \text{otherwise},
\end{array}\right.\eqe
where $\mu$ is a probability measure on $\N_0$, $q_0=\sum_{y\in\N}q_{0,y}$, and $r(x)$ is a positive finite function on $\N_0$ \cite{C97} .
Assume

\medskip
\noindent($\rm\mathbf{H1}$) $\mu(0)>0$, $\mu(0)+\mu(1)<1$.

\medskip
\noindent($\rm\mathbf{H2}$) $\sum_{y\in\N}q_{0,y}<\infty$.
\medskip

\thmb
Assume ${\rm(\mathbf{H1})}$-${\rm(\mathbf{H2})}$. Let $(Y_t\colon t\ge 0)$ be a process generated by the $Q$-matrix given in \eqref{Eq-5} and $Y_0\neq0$. Then $(Y_t\colon t\ge 0)$ is irreducible if $q_0>0$ and the conditional process of $(Y_t\colon t\ge 0)$ before absorption is irreducible if $q_0=0$.
\thme
\prb
First, assume $q_0=0$. ${\rm(\mathbf{H1})}$-${\rm(\mathbf{H2})}$ imply ${\rm(\mathbf{A1})}$-${\rm(\mathbf{A2})}$ are satisfied. Since $\cT_-=\{-1\}$ and $\cT_+=\supp\mu\setminus\{0,1\}-1$, we have $\om^*=1$. Hence the irreducibility of the conditional process follows from Corollary~\ref{co-1}.

Assume $q_0>0$, then $0$ communicates with states in $\N$, and thus $Y_t$ is irreducible by the above analysis for $q_0=0$.
\pre

\section{Proof of Theorem \ref{th-1}}

None of the jumps in $\cT$ are active in a state in $\N_0^d\setminus\overline{\cI}$, i.e., $\N_0^d\setminus\overline{\cI}\subseteq\sT\cup\sN$. On the other hand, if $x\in \sT\cup\sN$, then $\lambda_{\om}(x)=0$ for $\om\in\cT$, hence $\sT\cup\sN\subseteq\N_0^d\setminus\overline{\cI}$, and  equality holds. Hence $\sE\cup\sP\cup\sQ=\overline{\cI}$.  It now suffices to show that $\sN=\N_0^d\setminus\lt(\overline{\cO\cup\cI}\rt)$, based on the basic property $\overline{A}\cup\overline{B}=\overline{A\cup B}$ and $\overline{A}\setminus \overline{B}=\overline{A\setminus B}\setminus \overline{B}$. First, it is obvious that  $\sN\subseteq\N_0^d\setminus\overline{\cO\cup\cI}$. Conversely, suppose there exists $x\in\N_0^d\setminus\overline{\cO\cup\cI}$ such that $x\rightharpoonup y$ for some $y\in\N_0^d$. {Then there must exist a path from $x$ to $y$, which implies that there exists a jump active in $x$, i.e., $x\in \overline{\cI}$, a contradiction. Analogously, one can show that none of the states in $\N_0^d\setminus\overline{\cO\cup\cI}$ are reachable from any state in $\N_0^d$. This shows that $\N_0^d\setminus\overline{\cO\cup\cI}\subseteq\sN$.
Next, we verify the expressions for $\sE$ and $\sP\cup\sQ$ separately.
By the definition of escaping states, $\sE=\{x\in\overline{\cI}\colon x\rightharpoonup y\ \text{implies}\ y\not\rightharpoonup x\}$, and for any $x\in\overline{\cI}_{\om}^o$, $x \xldownrupharpoon[]{}  x+\om$, we have $\sE\subseteq\overline{\cI}\setminus\bigcup_{\om\in\cT
}\overline{\cI}_{\om}^o$. Conversely, for any $x\in\overline{\cI}\setminus\bigcup_{\om\in\cT
}\overline{\cI}_{\om}^o$, suppose $x \xldownrupharpoon[]{} y$ for some $y\in\N^d_0$. Then there exists a cycle through both $x$ and $y$, and thus there exists $\om\in\cT$ such that $x \xldownrupharpoon[]{} x+\om$, i.e., $x\in\overline{\cI}_{\om}^o$, a contradiction.
Since $\sE\cup\sP\cup\sQ=\overline{\cI}$ and $\bigcup_{\om\in\cT}\overline{\cI}_{\om}^o\subseteq\overline{\cI}$, we have
$\sP\cup\sQ=\bigcup_{\om\in\cT}\overline{\cI}_{\om}^o$.

For any $x\in\overline{\cI}_{\om}\setminus\overline{\cI}_{\om}^o$ for some $\om$, $x\rightharpoonup x+\om$ but $x+\om\not\rightharpoonup x$. Hence $x\in\sE\cup\sQ$, within an open communicating class, i.e., $\bigcup_{\om\in\cT\setminus\cT^o}\bigl(\overline{\cI}_{\om}\setminus\overline{\cI}_{\om}^o\bigr)\subseteq\sE\cup\sQ$. Finally, we show the remaining conclusions one by one.

(i) We first prove the sufficiency. It suffices to show $\sE\cup\sQ=\varnothing$, which further implies that $\sT=\varnothing$, and thus $\sP=\overline{\cI}$ and $\sN=\N^d_0\setminus\overline{\cI}$. Suppose $\sE\cup\sQ\neq\varnothing$. Then there exist $x, y\in\overline{\cI}$ such that
$x\rightharpoonup y,\ y\not\rightharpoonup x$.
Assume
\[x=x^{(1)}\rightharpoonup_{\widetilde{\om}^{(1)}}x^{(2)}\rightharpoonup_{\widetilde{\om}^{(2)}}\cdots\rightharpoonup_{\widetilde{\om}^{(\widetilde{m}-2)}} x^{(\widetilde{m}-1)}\rightharpoonup_{\widetilde{\om}^{(\widetilde{m}-1)}}x^{(\widetilde{m})}=y.\]
By the condition,  $x^{(j)} \xldownrupharpoon[]{} x^{(j+1)}$ for all $j=1,\cdots,\widetilde{m}-1$. This further implies that $x \xldownrupharpoon[]{} y$, a contradiction. Next we prove the necessity. Since ${\bigcup}_{\om\in\cT}\bigl(\overline{\cI}_{\om}\setminus\overline{\cI}_{\om}^o\bigr)\subseteq\sE\cup\sQ=\varnothing$, it follows that $\cT^o=\cT$.

(ii) It  follows from the expression for $\sP\cup\sQ$ as well as the definition of $\cI_{\om}^o$.

(iii) It suffices to show $\overline{\cI}_{\om}^o\subseteq\overline{\cO}$ for all $\om\in\cT$. Let $x\in\overline{\cI}_{\om}^o$, then $x+\om\rightharpoonup x$, and thus we have $x\in\overline{\cO}_{\widetilde{\om}}$ for some $\widetilde{\om}$, i.e., $\overline{\cI}_{\om}^o\subseteq\overline{\cO}$.

(iv) Based on the above analysis, $\N^d_0\setminus\bigcup_{\om\in\cT
}\overline{\cI}_{\om}^o=\sN\cup\sT\cup\sE$, consists of all singleton communicating classes, and $\sP\cup\sQ=\bigcup_{\om\in\cT}\overline{\cI}_{\om}^o$. If the condition holds, then every $x\in\sP\cup\sQ$ leads to a singleton communicating class. Hence $x$ must be within an open communicating class, and thus $x\notin\sP$. This shows that $\sP=\varnothing$, and $\sQ=\bigcup_{\om\in\cT}\overline{\cI}_{\om}^o$.

(v) Since $\bigcup_{\om\in\cT}\overline{\cI}_{\om}^o=\sP\cup\sQ\neq\varnothing$, it contains at least two states. The conclusion follows from the definition of reachability.
\smallskip

To show that all communicating classes are finite, it suffices to show that the underlying invariant subspaces defined in \eqref{Eq-subspace} are all finite, since every communicating class lies in one subspace. Let $c\in\N^d_0$. By the assumption, $\nu\cdot\om=0$ for all $\om\in\cT$. For all $x\in\sL_c$, $x-c\in\sS$, and thus $(x-c)\cdot\nu=0$. Hence $\sum_{j=1}^dx_j\nu_j=\sum_{j=1}^dc_j\nu_j$. Since $\nu\in\N^d$, we have $x_i\le\frac{\sum_{j=1}^dc_j\nu_j}{\nu_i}$, for all $i=1,\ldots,d$, and hence $\sL_c$ is finite.

\section{Proof of Theorem \ref{th-6}}

The proof relies on four lemmata provided in the appendix.
We only prove Theorem \ref{th-6} for the case $\om^*\in\Z^d\setminus\N^d_0.$  The proof can readily be adapted to the  case $\om^*\in\N^d_0$ with $b=\bz$. Indeed, for all $c\in\N_0^d$, by translational invariance, one can find $c'$ large enough in all coordinates $j$ for $\om^*_j>0$ such that $\sL_c=\sL_{c'}$. Then it  suffices to replace $c$ by $c'$ in the rest of the proof.

By Lemma~\ref{le-1}, choose a finite $\widetilde{\cT}\subseteq\cT$ such that
$$\underset{\om\in\widetilde{\cT}_+}{\cup}\overline{\cI}_{\om}=\overline{\cI}_+,\ \underset{\om\in\widetilde{\cT}_-}{\cup}\overline{\cI}_{\om}=\overline{\cI}_-,\ \underset{\om\in\widetilde{\cT}}{\cup}\overline{\cI}_{\om}=\overline{\cI},\ \underset{\om\in\widetilde{\cT}_+}{\cup}\overline{\cO}_{\om}=\overline{\cO}_+,\ \underset{\om\in\widetilde{\cT}_-}{\cup}\overline{\cO}_{\om}=\overline{\cO}_-,\ \underset{\om\in\widetilde{\cT}}{\cup}\overline{\cO}_{\om}=\overline{\cO},$$
where $\widetilde{\cT}_{\pm}=\{\omega\in\widetilde{\cT}\colon {\rm sgn}(\om_1)=\pm1\}$.  Let $\widetilde{\cI}=\cup_{\om\in\widetilde{\cT}}\cI_{\om}$ and $\widetilde{\cO}=\cup_{\om\in\widetilde{\cT}}\cO_{\om}$.

Let $M=\underset{\om\in\widetilde{\cT}}{\max}\lt|\frac{\om}{\om^*}\rt|+1$ ($\ge 2$) and define $b\in\N_0^d$ by
  $$b_j=M|\om^*_j|+{\max}(\widetilde{\cI}_j\cup\widetilde{\cO}_j),\quad \text{for}\quad j=1,\ldots,d.$$

   Hence\eqb\label{Eq-10}\om^*[-M,M]+b\subseteq\underset{y\in\widetilde{\cI}\cup\widetilde{\cO}}{\cap}\overline{\{y\}},\eqe which implies that all jumps in $\widetilde{\cT}$ are active in all states in $\om^*[-M,M]_1+b$. Let $c\in \N_0^d+b$. For convenience, within this proof we slightly abuse $\overline{\cI}$ to mean $\overline{\cI}\cap\sL_c$ to ignore the dependence on $c$. Analogously for $\overline{\cI}_+\cap\sL_c$, etc. Let
  $$D_k=\om^*[-M+1,M-1]_1+\frac{k-1}{\om^{**}}\om^*+c,\quad \text{for}\quad k=1,\ldots,\om^{**}.$$
    Given $k\in[1,\om^{**}]_1$. By the definition of $\widetilde{\cT}$, there exist $-m_1\om^*\in\widetilde{\cT}_-$ and $m_2\om^*\in\widetilde{\cT}_+$ with $m_1,\ m_2\in\N$ coprime. Since $c\ge b$, all jumps in $\widetilde{\cT}$ including $-m_1\om^*$ and $m_2\om^*$ are active in every state in $D_k$, in the light of \eqref{Eq-10}. By Lemma~\ref{Sle-0}, $D_k$ is communicable.

Moreover, by Lemma \ref{Sle-8}, the sets $\overline{\cI}_{+},\  \overline{\cI}_{-},\ \overline{\cO}_{+}$ and $\overline{\cO}_{-}$ are all non-empty lattice intervals, and so are their finite intersections, and as well as their finite unions due to the non-emptiness of intersections.
In particular, $\sK_c=(\overline{\cI}_+\cap \overline{\cO}_-)\cup(\overline{\cI}_-\cap \overline{\cO}_+)$ is also a lattice interval.
Hence $\sK_c=[\uc,\bc]_1$.

For every $k\in[1,\om^{**}]_1$, let $G_k=\sK_c\cap \Gamma^{(k)}_c$. In  Step I and Step II below, we will show that $G_k$ is a communicating class with at least two distinct states, which in turn implies that $G_k$ is either a PIC or a QIC. In Step III, we show that $G_k$ is a QIC trapped into $\sT_c^{(k)}$ for all $k\in\Sigma^+_c$, and a PIC for all $k\in \Sigma^-_c$.

\medskip
\noindent {\bf Step I}. $G_k$ is communicable with $\# G_k\ge2$.
Indeed, since $M\ge2$ and $c\ge b$,
$$\lt(c\pm(M-1)\om^*+\frac{k-1}{\om^{**}}\om^*\rt)_j\ge {\max}(\cI_j\cup\cO_j),  \quad \forall j=1,\ldots,d,$$
which implies that $$c\pm(M-1)\om^*+\frac{k-1}{\om^{**}}\om^*\in\Gamma_c^{(k)}\cap(\overline{\cI}_+\cap \overline{\cO}_-)\cap(\overline{\cI}_-\cap \overline{\cO}_+)\subseteq G_k.$$
This shows  $\# G_k\ge2$.

Next we prove that $G_k$ is communicable. Let $D=\cup_{l=1}^{\om^{**}}D_l$. By Lemma \ref{Sle-8}, $D=[\uD,\bD]_1$ is a lattice interval with $\uD=c+(-M+1)\om^*$ and $\bD=c+(M-1)\om^*+\frac{\om^{**}-1}{\om^{**}}\om^*$. By Lemma \ref{Sle-6}, for any $x\in D_i$ and any $y\in D_j$ with $i\neq j$, $i, j\in[1,\om^{**}]_1$, $x$ neither is reachable from nor leads to $y$. Since both $D$ and $\sK_c$ are lattice intervals, we have $\sK_c\setminus D=[\uc,\uD[_1\cup]\bD,\bc]_1$. Due to the communicability of $D_k$, to see $G_k$ is communicable, it suffices to show for all $x\in \sK_c\setminus D$, there exists $y\in D$ such that $x \xldownrupharpoon[]{} y$.

In the following, we prove that for all $x\in[\uc,\uD[_1$, there exists $y\in D$ such that $x \xldownrupharpoon[]{} y$. The analogous property holds for $]\bc,\bD]_1$.

Recall that $\sK_c=(\overline{\cI}_+\cap \overline{\cO}_-)\cup(\overline{\cI}_-\cap \overline{\cO}_+)$, and both $\overline{\cI}_+\cap \overline{\cO}_-$ and $\overline{\cI}_-\cap \overline{\cO}_+$ are lattice intervals.
Since $D\subseteq\overline{\cI}_+\cap \overline{\cO}_-\cap\overline{\cI}_-\cap \overline{\cO}_+$, assume w.o.l.g. that $\varnothing\neq[\uc,\uD[_1\subseteq\overline{\cI}_+\cap \overline{\cO}_-$. Then $\uc=\min_1\overline{\cI}_+\cap \overline{\cO}_-$ and $[\uc,\bD]_1\subseteq\overline{\cI}_+\cap \overline{\cO}_-$.
Note that
$$\tfrac{\bD-\uD}{\om^*}>M-1=\underset{\om\in\widetilde{\cT}}{\max}\lt|\tfrac{\om}{\om^*}\rt|,$$
 and thus every state in $[\uc,\uD[_1$ cannot lead to a state in  $]\bD,\bc]_1$ within a single jump. Let $x\in[\uc,\uD[_1$. On the one hand, since $x\in[\uc,\bD]_1\subseteq\overline{\cI}_+$, and one can show by induction that there exists $y\in D$ such that $x\rightharpoonup y$, realized by a finite ordered set of jumps in $\cT_+$, with $y-x\in\om^*\N$. On the other hand, since $x\in \overline{\cO}_-$, in an analogous manner, one can show that there exists $z\in D$ such that $z\rightharpoonup x$, realized by a finite ordered set of jumps in $\cT_-$, with $-(x-z)\in\om^*\N$. Hence $y-z=y-x+x-z\in\om^*\Z$. By Lemma~\ref{Sle-6}, $y, z\in D_k$ for some $k\in[1,\om^{**}]_1$, i.e., $y=z$ or $y \xldownrupharpoon[]{} z$. By transitivity, $x \xldownrupharpoon[]{} y\in D$.

\medskip
\noindent {\bf Step II}. $\overline{\cI}\setminus \sK_c=\sE$. From Theorem~\ref{th-1}, $\sT_c\cup\sN_c=\sL_c\setminus \overline{\cI}$, and thus it suffices to show that $\sL_c\setminus \sK_c$ is composed of singleton communicating classes assuming that $\overline{\cI}\setminus \sK_c\neq\varnothing$.

Since $\sK_c\subseteq\overline{\cI}$, and $\sK_c$ and $\overline{\cI}$ are both lattice intervals, we have $\overline{\cI}\setminus \sK_c=[\min_1\overline{\cI},\uc[_1\cup\, ]\bc,\max_1\overline{\cI}]_1$. Assume w.o.l.g. that $[\min_1\overline{\cI},\uc[_1\neq\varnothing$. It then suffices to show that $[\min_1\overline{\cI},\uc[_1$ is composed of singleton communicating classes. It is easy to see that
$$\min\nolimits_1\overline{\cI}_+<_1\min\nolimits_1\overline{\cO}_+,\quad \min\nolimits_1\overline{\cI}_->_1\min\nolimits_1\overline{\cO}_-.$$
Since $\overline{\cI}_{+},\  \overline{\cI}_{-},\ \overline{\cO}_{+}$ and $\overline{\cO}_{-}$ are all lattice intervals,  it readily yields that
$$ [\min\nolimits_1\overline{\cI},\uc[_1\subseteq\overline{\cI}_-\setminus(\overline{\cI}_+\cup\overline{\cO}_+)\ \text{or}\ [\min\nolimits_1\overline{\cI},\uc[_1\subseteq\overline{\cI}_+\setminus(\overline{\cI}_-\cup\overline{\cO}_-).$$
Further assume w.o.l.g. that $[\min_1\overline{\cI},\uc[_1\subseteq\overline{\cI}_-\setminus(\overline{\cI}_+\cup\overline{\cO}_+)$. Let $x\in[\min_1\overline{\cI},\uc[_1$. Now we only show that no other state communicates with $x$ by contradiction. Suppose there exists $y\neq x$ such that $x \xldownrupharpoon[]{} y$. Then there exists a cycle connecting $x$ and $y$, denoted by
$$x=y^{(0)}\rightharpoonup\ldots\rightharpoonup y^{(m)}\rightharpoonup y^{(0)}.$$
Let $z={\min_1} \{y^{(j)}: 0\le j\le m\}=y^{(k)}$ for some $0\le k\le m$. Since  $k\not\equiv (k+1)\hskip-0.2cm\mod\hskip-0.1cm(m+1)$, we have $z<_1y^{(k+1)\hskip-0.2cm\mod\hskip-0.1cm(m+1)}$, and thus $z\in \overline{\cI}_+$, for $z\rightharpoonup y^{(k+1)\hskip-0.2cm\mod\hskip-0.1cm(m+1)}$ must be realized by a jump in $\cT_+$. On the one hand, since $x\notin \overline{\cI}_+$, we have  $z\le_1 x<_1\min_1\overline{\cI}_+$, and thus $z\notin \overline{\cI}_+$, recalling again that $\overline{\cI}_+$ is a lattice interval. This is a contradiction.

\medskip
\noindent {\bf Step III}. Given $k\in\Sigma^+_c$, for every $x\in \sQ^{(k)}_c$ and $y\in\sT^{(k)}_c$, we have $x\rightharpoonup y$.

Based on Steps I and II, by Lemma \ref{Sle-6}, $\sQ^{(k)}_c$ is a
quasi-irreducible component for all $k\in\Sigma^+_c$, and $\sP^{(k)}_c$ a positive irreducible component for all $k\in\Sigma^-_c$. In particular, there are precisely $\# \Sigma^+_c$ quasi-irreducible components ultimately leading only to trapping states, and there are $\# \Sigma_c^-$ positive irreducible components.

In the light of Lemma \ref{Sle-6}, it suffices to show that:  $\forall x\in\sT$, there exists $y\in \sK_c$ such that $y\rightharpoonup x$.

Again on account of Lemma \ref{Sle-6}, we assume w.o.l.g. that $\sT\neq\varnothing$ and $\om^{**}=1$.
Given $x\in\sT$, there exists $z\in\cI=\sE\cup \sK_c$ such that $z\rightharpoonup x$. Assume w.o.l.g. that $z\in\sE$.
From Step II,
$$z\in\overline{\cI}\setminus \sK_c=[\min\nolimits_1\overline{\cI},\uc[_1\cup]\bc,\max\nolimits_1\overline{\cI}]_1.$$
Furthermore assume w.o.l.g. that $z\in[\min_1\overline{\cI},\uc[_1$. By the analysis in Step II, it suffices to prove that there exists $y\in \sK_c$ such that $y\rightharpoonup z$ under the further assumption w.o.l.g. that
$[\min_1\overline{\cI},\uc[_1\subseteq\overline{\cI}_-\setminus\overline{\cI_+\cup\cO_+}.$
 Since $\min_1\overline{\cI}_->_1\min_1\overline{\cO}_-$, we have $z\in\overline{\cO}_-$, and
 $$\tfrac{\bc-\uc}{\om^*}>\underset{\om\in\widetilde{\cT}}{\max}\lt|\tfrac{\om}{\om^*}\rt|,$$
  similarly to Step I, one can show by induction that there exists $y\in \sK_c$ such that $y\rightharpoonup z$, realized by a finite ordered set of jumps in $\widetilde{\cT}_-$.

\section{Proof of Corollary~\ref{co-1}}

We assume w.o.l.g. that $d=1$ and $c=0$. Therefore $\om^*=\om^{**}$, and the partial order indexed by `1' coincides with the natural partial order on the real line. Hence the dependence on $1$ and $c$ is omitted in the notation of the sets, numbers as well as symbols (e.g., in $\Gamma_c^{(k)}$, $\so(c)$, and $\ge_1$, etc.). Moreover, $c^*=\infty$, $\min_1\sL_c=c=0$, and $\Gamma^{(k)}=\om^*\N_0+k-1$, for $k=1,\ldots,\om^*$.

Let $\so_+=\min\overline{\cO}_+$ and $\si_-=\min\overline{\cI}_-$. Hence $$\overline{\cI}=\N_0+\si,\ \overline{\cI}_+=\N_0+\si_+,\ \overline{\cI}_-=\N_0+\si_-,\ \overline{\cO}=\N_0+\so,\  \overline{\cO}_+=\N_0+\so_+,\ \overline{\cO}_-=\N_0+\so_-.$$
Moreover, $\si=\min\{\si_-,\si_+\},\ \so=\{\so_-,\so_+\},\ \si_+<\so_+,\ \si_->\so_-$. When $\so\le\si$, then $\so\le\si\le\si_+<\so_+$, which implies that $\so=\so_-<\si_-$. Hence  $\min\{\si,\so\}=\so_-=\min\{\si_+,\so_-\}$. When $\so\ge\si$, then $\si\le\so\le\so_-<\si_-$, which implies that $\si=\si_+$ and hence $\min\{\si,\so\}=\si_+=\min\{\si_+,\so_-\}$. Hence it always holds that
$\min\{\si,\so\}=\min\{\si_+,\so_-\}$.

The expressions for $\sN$, $\sT$ and $\sE$ follow immediately from Theorem~\ref{th-6}. Indeed, since $\overline{\cO}\cup\overline{\cI}=\N+\min\{\si,\so\}$, $\sN=\N_0\setminus(\overline{\cO}\cup\overline{\cI})=[0,\ldots,\min\{\si,\so\}-1[_1$.
Moreover, , and $\sK_c=\N_0+c_*$, $\sT=\overline{\cO}\setminus\overline{\cI}=(\N_0+\so)\setminus(\N_0+\si)=[\so,\ldots,\si[_1$,
and $\sE=\overline{\cI}\setminus\sK_c=[\si,c_*[_1$, where
\[\begin{split}
c_*=&\min\{\min(\cI_+\cap\cO_-),\min(\cI_-\cap\cO_+)\}\\
=&\min\{\max\{\si_+,\so_-\},\max\{\si_-,\so_+\}\}\\
=&\max\{\si_+,\so_-\},
\end{split}\]
as $\si_+\le\so_+$ and $\so_-\le\si_-$.
This verifies that expression for $\sE$.
Moreover, also  from Theorem~\ref{th-6}, it follows that $\sP\cup\sQ=\sK_c=\N_0+c_*=\N_0+\max\{\si_+,\so_-\}$.

Next, we express each PIC and QIC explicitly.

\noindent\rm{(1)} Assume $\sT\neq\varnothing$. Then $\so<\si\le\si_+<\so_+$, and hence $\so=\so_-$.
In this case, $\so_-<\si_+$, and therefore $\min\{\si_+,\so_-\}=\so_-$, $\max\{\si_+,\so_-\}=\si_+$, and
\[\sN=\{0,\ldots,\so-1\},\quad \sT=\{\so,\ldots,\si-1\},\quad \sE=\{\si,\ldots,\si_+-1\},\quad \sP\cup\sQ=\N_0+\si_+.\]
Let $\widetilde{\Sigma}^+=\{k\in\{1,\ldots,\om^*\}\colon \sT^{(k)}\neq\varnothing\}$. In the following, we will show $\widetilde{\Sigma}^+=\Sigma^+$.
Let $k\in\widetilde{\Sigma}^+$ and $x\in\sT^{(k)}$. Then $x\ge\so_-$. By ($\rm\mathbf{A1}$), there exists $y\in\N_0+\si_+$ such that $y\rightharpoonup x$. Hence $y\in\sQ$ is in an open class. This further shows $\sQ^{(k)}\neq\varnothing$, and hence $\widetilde{\Sigma}^+\subseteq\Sigma^+$. On the other hand, for any $k\in\Sigma^+\setminus\widetilde{\Sigma}^+$, $\sT^{(k)}=\varnothing$ while $\sQ^{(k)}\neq\varnothing$, this means $\sP^{(k)}=\varnothing$, i.e., there is no closed class in $\Gamma^{(k)}$, which is impossible since $\Gamma^{(k)}$ is closed by Lemma~\ref{Sle-6}. Hence $\Sigma^+=\widetilde{\Sigma}^+$ (this holds trivially with the same argument when $\sT=\varnothing$).

Next we show the expressions for $\Gamma^{(k)}\cap(\om^*\N_0+c_*)$, and determine the sets $\Sigma^+$ and $\Sigma^-$. Let $k=1,\ldots,\om^*$.
On one hand, for any  $x\in\om^*\N_0+\lc(\si_+-k+1)/\om^*\rc\om^*+k-1$, we have $$x\ge\lc(\si_+-k+1)/\om^*\rc\om^*+k-1\ge\si_+=\max\{\si_+,\so_-\}=c_*,$$ which implies that
 \[\begin{split}
\om^*\N_0+\lc(\si_+-k+1)/\om^*\rc\om^*+k-1
\subseteq\Gamma^{(k)}\cap(\om^*\N_0+c_*).
\end{split}\]
On the other hand, for any $x\in\Gamma^{(k)}\cap\om^*\N_0+c_*$, we have $x\in\om^*\N_0+k-1$ and $x\ge\si_+$. Hence $x-k+1\ge\si_+-k+1$ and $(x-k+1)/\om^*\in\N_0$, which implies that $(x-k+1)/\om^*\ge\lc(\si_+-k+1)/\om^*\rc$, i.e., $x\in\om^*\N_0+\lc(\si_+-k+1)/\om^*\rc\om^*+k-1$. Hence we show $\om^*\N_0+\lc(\si_+-k+1)/\om^*\rc\om^*+k-1=\Gamma^{(k)}\cap\om^*\N_0+c_*$. For every $v\in[\so,\ldots,\max\{\si,\so+\om^*\}[_1$, let $k_v=1+v-\lf v/\om^*\rf\om^*$. Hence $1\le k_v\le\om^*$ and $1+v-k_v\in\om^*\N_0$. This shows that $k_v\in\Sigma^+$, i.e., $\sT^{(k_v)}\neq\varnothing$, if and only if  $v\in\sT$. Recall that $\sT=[\so,\si[_1$. Therefore $$\Sigma^+=\{1+v-\lf v/\om^*\rf\om^*\colon v\in[\so,\si[_1\},\quad \Sigma^-=\{1+v-\lf v/\om^*\rf\om^*\colon v\in[\si,\so+\om^*[_1\},$$ and the expressions for $\sP^{(k)}$ and $\sQ^{(k)}$ follow directly by definition.

\noindent{\rm(2)} Assume $\sT=\varnothing$. Then $\so\ge\si$, and hence $\si_->\so_-\ge\so\ge\si$, implying $\si=\si_+\le\so\le\so_-$. In this case, $\max\{\si_+,\so_-\}=\so_-$ and $\min\{\si_+,\so_-\}=\si_+$, and $$\sN=\{0,\ldots,\si_+-1\},\quad \sE=\{\si_+,\ldots,\so_--1\},\quad \sP\cup\sQ=\N_0+\so_-.$$ Moreover, $\Sigma^+=\varnothing$ and $\Sigma^-=\{1,\ldots,\om^*\}$.  The rest argument to verify the expressions for $\sP^{(k)}$ is the same as (1).

\medskip

We continue to prove the rest statements.
(i) When $\#\sT\ge\om^*$, i.e., $\si-\so\ge\om^*$, we have $\Sigma^+=[1,\om^*]_1$, and hence $\#\Sigma^+=\om^*$.
(ii) When $0<\#\sT<\om^*$, i.e., $0<\si-\so<\om^*$, we have $\#\Sigma^+=\#\sT=\si-\so$.
(iii) When $\#\sT=0$, i.e., $\sT=\varnothing$, we have $\si\le\so$, and $\Sigma^+=\varnothing$ and $\#\Sigma^+=0$.
In all cases, $\# \Sigma^+={\min}\lt\{\om^{*},\max\{0,\si(c)-\so(c)\}\rt\}$.

Since $\widetilde{\Sigma}^+=\Sigma^+$, we have $\Sigma^+\neq\varnothing$ implies $\sT\neq\varnothing$. Moreover, in this case, $\max\sE<\si_+$. By the definition of $\si_+$, no escaping state can reach any state in $\sQ\cup\sP$,  but only can reach another escaping state or a trapping state, since $\sE$ is open.

\section*{Acknowledgements}

The work was initiated with most part of it done when the first author was at the University of Copenhagen. The authors thank the editors' and referees' comments which helped improve the presentation of the paper. The authors acknowledge the support from The Erwin Schr\"{o}dinger Institute (ESI) for the workshop on ``Advances in Chemical Reaction Network Theory''. CX acknowledges the support from TUM University Foundation and the Alexander von Humboldt Foundation. CW acknowledges support from the Novo Nordisk Foundation (Denmark), grant NNF19OC0058354.

\appendix

\renewcommand\appendixname{Appendix}

\section{Basic facts in set theory}

To be self-contained, we provide two elementary results in set theory.

\prob\label{Spro-8}
Let $A\subseteq\Z^d\setminus\{\bz\}$ be non-empty. Then $\dim\spa A=1$ if and only if $A$ has a common divisor if and only if $\gcd(A)$ exists.
\proe
\prb
If $\gcd(A)$ exists,  then clearly $A$ has a common divisor. If $A$ has a common divisor $a$, then necessarily $\spa A\subseteq a\Z$. Hence, $\dim\spa A\le1$. On the other hand, $A\neq\{\bz\}$ and $A$ is non-empty, hence $\dim\spa A\ge1$. So $\dim\spa A=1$. If $\dim\spa A=1$, then there exists $a\in\Z^d\setminus\{\bz\}$ such that $\spa A=a\Z$. Clearly, $a|b$, $\forall b\in\spa A$, in particular $a|b$, $\forall b\in A$. Hence, $A$ has a common divisor $a$. Since $A\neq\{\bz\}$ and $A$ is non-empty, then the set of common divisors is finite, and denoted by $\{a_1,\ldots,a_k\}$. We will show by contradiction that $A$ has a gcd. Assume $a_i\le a_k$ (recall $\le$ respects the order by the first coordinate) for all $1\le i\le k$ and that $A$ does not have a gcd. Then $\tfrac{a_k}{a_i}=\frac{q}{p}$, $p,q\in\Z$, for some $i$, such that the fraction is irreducible. Since $a_i|b$ and $a_k|b$ for $b\in A$, then it must be that also $pa_k|b$, contradicting that $a_k$ is the largest among all common divisors . Hence $A$ has a gcd.
\pre

\begin{proposition}
\label{Spro-11}
Let $B$ be a non-empty subset of $\N^d_0$. Then $B$ has a finite non-empty minimal set $E$, and $B\subseteq \overline{E}$.
\end{proposition}
\prb
We first prove the former part of the conclusion.
That $E\neq\varnothing$ follows directly from Zorn's lemma. In the following we show $\#E<\infty$. Suppose $\#E=\infty$. Then there exists $1\le k\le d$ and a sequence $(x^{(j)})_{j\ge1}$ such that  $(x^{(j)}_k)_{j\ge1}$ is unbounded. Assume w.o.l.g. $k=1$ and $x^{(j)}_1\uparrow\infty$ as $j\to\infty$. If the remaining sets  $(x^{(j)}_i)_{j\ge1}$, for $i=2,\ldots,d$, are bounded, say $x^{(j)}_i\le M$ for $j\ge 1$, then there will be two points such that $x^{(j_1)}\ge x^{(j_2)}$, contradicting minimality. Hence,  we might assume w.l.o.g. that $x^{(j)}_i\uparrow\infty$, for $i=1,2$, as $j\to\infty$. Continuing in this fashion yields a sequence $(x^{(j)})_{j\ge1}$ such that $x^{(j)}\le x^{(j+1)}$, $j\ge 1$, contradicting minimality. Hence, $E$ is finite.

Next we prove $B\subseteq \overline{E}$. Suppose there exists $x^{(1)}\in B\setminus \overline{E}\subseteq B\setminus E$. By minimality, there exists $x^{(2)}\in B\setminus\{x^{(1)}\}$ such that $x^{(1)}\ge x^{(2)}$. Since $x^{(1)}\notin \overline{E}$, we have $x^{(2)}\in B\setminus \overline{E}$. By induction one can get a decreasing sequence $(x^{(j)})_{j\ge 1}\subseteq B\setminus \overline{E}$ of distinct elements. This is impossible since  $\#(x^{(j)})_{j\ge 1}\le\prod_{k=1}^d(x^{(1)}_k+1)$.
\pre

\section{Lemmata for Theorem~\ref{th-6}}

Recall that  $\dim\sS=1$, $\rm(\mathbf{A1})$-$(\rm\mathbf{A2})$ are assumed for these lemmata, and that $\om^*=\gcd(\cT)$.
The first lemma establishes a result on communicability of two states.

\leb\label{Sle-0}
Assume $\om^{(1)}=-m_1\om^*\in\cT_-$ and $\om^{(2)}=m_2\om^*\in\cT_+$ with coprime $m_1,\ m_2\in\N$. Let $x\in\N^d_0$. If $\om^{(1)}$ is active in $x+(m_1-1+j)\om^*$ for all $j\in[1,m_2]_1$ and $\om^{(2)}$ is active in $x+(j-1)\om^*$ for all $j\in[1,m_1]_1$, then $\om^*[0,m_1+m_2-1]_1+x$ is communicable. In particular, for all integers $M\ge m_1+m_2-1$, if $\om^{(1)}$ and $\om^{(2)}$ are both active in $\om^*[0,M]_1$, then $\om^*[0,M]_1+x$ is communicable.
\lee

\prb
Since the following arguments are independent of the specific form of $\om^*$, assume w.o.l.g. that $d=1$, $\om^*=1$ and $x=0$. It then suffices to show that $[0,m_1+m_2-1]_1$ is communicable. We first show $[0,\max\{m_1,m_2\}-1]_1$ is communicable, and then show $[0,m_1+m_2-1]_1$ is communicable.

\medskip
\noindent {\bf Step I.}
$[0,\max\{m_1,m_2\}-1]_1$ is communicable. To see this, we first show that  $[0,m_2-1]_1$ is communicable. Then in an analogous way, one can show $[0,m_1-1]_1$ is also communicable. This implies that $[0,\max\{m_1,m_2\}-1]_1$ is communicable.

Indeed, there exists a cycle connecting all states in $[0,m_2-1]_1$, which immediately yields that $[0,m_2-1]_1$ is communicable.
To prove this, one needs the following elementary identity,
\begin{equation}\label{SEq-11}
\lt\{\lt\lceil\frac{jm_1}{m_2}\rt\rceil m_2-jm_1\colon j\in\N\rt\}=\lt\{\lt\lceil\frac{jm_1}{m_2}\rt\rceil m_2-jm_1\colon j\in [1,m_2]_1\rt\}=[0,m_2-1]_1,
\end{equation}
(a proof is included below)
from which it immediately follows that
\[0\le\lt\lceil\frac{jm_1}{m_2}\rt\rceil m_2-jm_1\le \lt\lceil\frac{jm_1}{m_2}\rt\rceil m_2-(j-1)m_1\le m_1+m_2-1,\]
which further yields the desired cycle  by repeated jumps of $\om^{(1)}$ and $\om^{(2)}$, connecting the states in
$\lt\{\lt\lceil\frac{jm_1}{m_2}\rt\rceil m_2-jm_1\colon j\in [1,m_2]_1\rt\}$,

\[
\begin{split}0&\rightharpoonup \lt\lceil\frac{m_1}{m_2}\rt\rceil m_2\rightharpoonup\lt\lceil\frac{m_1}{m_2}\rt\rceil m_2-m_1\rightharpoonup \lt\lceil\frac{2m_1}{m_2}\rt\rceil m_2-m_1\rightharpoonup \lt\lceil\frac{2m_1}{m_2}\rt\rceil m_2-2m_1\\ &\rightharpoonup\cdots\rightharpoonup\lt\lceil\frac{jm_1}{m_2}\rt\rceil m_2-(j-1)m_1\rightharpoonup\lt\lceil\frac{jm_1}{m_2}\rt\rceil m_2-jm_1\\&\rightharpoonup\cdots\rightharpoonup\lt\lceil\frac{(m_2-1)m_1}{m_2}\rt\rceil m_2-(m_2-2)m_1\rightharpoonup\lt\lceil\frac{(m_2-1)m_1}{m_2}\rt\rceil m_2-(m_2-1)m_1\\
&\rightharpoonup\lt\lceil\frac{m_2m_1}{m_2}\rt\rceil m_2-(m_2-1)m_1\rightharpoonup \lt\lceil\frac{m_2m_1}{m_2}\rt\rceil m_2-m_2m_1\lt(=0\rt),\end{split}\] connecting all states $[0,m_2-1]_1$, which must be communicable.

For the reader's convenience, we here give a proof of the elementary identity \eqref{SEq-11}. On the one hand,
\[0\le\lt\lceil\frac{jm_1}{m_2}\rt\rceil m_2-jm_1<m_2,\quad \text{for all}\,\, j\in\N,\]
which yields
\eqb\label{SEq-12}\lt\{\lt\lceil\frac{jm_1}{m_2}\rt\rceil m_2-jm_1\colon j\in\N\rt\}\subseteq[0,m_2-1]_1.\eqe
On the other hand, for all $1\le i<j\le m_2$,
\eqb\label{SEq-13}\lt\lceil\frac{im_1}{m_2}\rt\rceil m_2-im_1\neq \lt\lceil\frac{jm_1}{m_2}\rt\rceil m_2-jm_1.\eqe
If  this is not so, then
 \[(j-i)m_1=\lt(\lt\lceil\frac{jm_1}{m_2}\rt\rceil-\lt\lceil\frac{im_1}{m_2}\rt\rceil\rt)m_2,\]
 for some $i<j$. Since $0<j-i<m_2$, and $m_1$, $m_2$ are coprime, $m_2\nmid(j-i)m_1$, a contradiction.

From \eqref{SEq-13}, it follows that
\[\#\lt\{\lt\lceil\frac{jm_1}{m_2}\rt\rceil m_2-jm_1\colon j\in [1,m_2]_1\rt\}\ge m_2,\]
which together with \eqref{SEq-12} yields \eqref{SEq-11}.

\medskip
\noindent

{\bf Step II.} $[0,m_1+m_2-1]_1$ is communicable,  since
 $$x\,\,\text{mod}\,\,m_1,\, y\,\,\text{mod}\,\,m_2\in [0,\max\{m_1,m_2\}-1]_1$$
 for $x,y\in [0,m_1+m_2-1]_1$,
and
 \[ x=\lt(x\,\,\text{mod}\,\,m_1\rt)+\left\lfloor\frac{x}{m_1}\right\rfloor m_1,\quad y=\lt(y\,\,\text{mod}\,\,m_2\rt)+\left\lfloor\frac{y}{m_2}\right\rfloor m_2. \]
Then,  similarly to Step I, a path from $x$ to $y$ is constructed.
The proof is complete.
\pre
The next lemma shows when two states are not reachable from each other.

\leb\label{Sle-6}
Let $c\in\N^d_0$. For any $x,y\in\sL_c$, if $x-y\notin\om^*\Z$, then $x\not\rightharpoonup y$ and $y\not\rightharpoonup x$.
\lee
\prb We prove $x\not\rightharpoonup y$ by contradiction. By symmetry,  $y\not\rightharpoonup x$. Suppose $x\rightharpoonup y$, then there exists a path from $x$ to $y$ realized by a finite ordered set of (possibly repeated) jumps $\{\om^{(j)}\}_{1\le j\le m}$ for some $m\in\N$ such that $y-x=\sum_{j=1}^m\om^{(j)}.$ Since $\om^{(j)}\in\om^*\Z$, then $y-x\in\om^*\Z$, which is a contradiction.
\pre

The following lemma guarantees connectedness of a discrete set.
\leb\label{Sle-8}
Let $U=\{u^{(j)}\}_{j=1}^{m}\subseteq\Z^d\setminus\{\bz\}$ for some $m\in\N$ and define $\bu$ by
$$\bu_i=\max_{1\le j\le m}u^{(j)}_i,\quad \forall 1\le i\le d.$$ Then for every $c\ge \bu$, $\sL_c\cap \overline{U}$ is a non-empty lattice interval.
\lee
\prb
Recall that if $A,\ B\subsetneq\sL_c$ are both lattice intervals, and $A\cap B\neq\varnothing$, then $A\cap B$ and $A\cup B$ are  also  lattice intervals. We use this property to prove the conclusion.
Indeed, \[\overline{U}=\overset{m}{\underset{j=1}{\bigcup}}\overline{\{u^{(j)}\}}.\]
If $c\ge\bu$, then $c\in \sL_c\cap \overline{\{u^{(j)}\}}$ for all $1\le j\le d$. Moreover, $\sL_c\cap \overline{\{u^{(j)}\}}$ is a lattice interval, since $\sL_c$ is a lattice interval, and $\overline{\{u^{(j)}\}}$ is a convex set in $\N^d_0$.
Hence $\sL_c\cap \overline{U}=\overset{m}{\underset{j=1}{\bigcup}}(\sL_c\cap \overline{\{u^{(j)}\}})$ is also a lattice interval containing $c$.
\pre
The following lemma guarantees that we can generate the sets $\cI$, etc., in terms of a finite set of jump vectors.
\leb\label{le-1}
Let $(\cT,\cF)$ be given with $\om^*=\gcd(\cT)$. Then there exists a finite set $\widetilde{\cT}\subseteq\Omega$ such that $\gcd(\widetilde{\cT})=\om^*$ and
$$\underset{\om\in\widetilde{\cT}_+}{\cup}\overline{\cI}_{\om}=\overline{\cI}_+,\ \underset{\om\in\widetilde{\cT}_-}{\cup}\overline{\cI}_{\om}=\overline{\cI}_-,\ \underset{\om\in\widetilde{\cT}}{\cup}\overline{\cI}_{\om}=\overline{\cI},\ \underset{\om\in\widetilde{\cT}_+}{\cup}\overline{\cO}_{\om}=\overline{\cO}_+,\ \underset{\om\in\widetilde{\cT}_-}{\cup}\overline{\cO}_{\om}=\overline{\cO}_-,\ \underset{\om\in\widetilde{\cT}}{\cup}\overline{\cO}_{\om}=\overline{\cO},$$
 where $\widetilde{\cT}_{\pm}=\{\omega\in\widetilde{\cT}\colon {\rm sgn}(\om_1)=\pm1\}$.
\lee

\prb
By Proposition~\ref{Spro-11}, $\overline{\cI}_+$ has a finite minimal set, say $E$, such that $\overline{\cI}_+\subseteq\overline{E}$. Since $\overline{\cI}_+=\cup_{\om\in\cT_+}\overline{\cI}_{\om}$, then for  $x\in E$, we have $x\in\overline{\cI}_{\om}$ for some $\om\in\cT_+$. By the definition of minimal element, $x\in\cI_{\om}$. Hence $E\subseteq\cup_{\om\in\cT'(\overline{\cI}_+)}\cI_{\om}$ for some finite set $\cT'(\overline{\cI}_+)\subseteq\cT_+$, implying that $$\underset{\om\in\cT'(\overline{\cI}_+)}{\cup}\overline{\cI}_{\om}\subseteq\cI_+\subseteq\overline{E}
\subseteq\underset{\om\in\cT'(\overline{\cI}_+)}{\cup}\overline{\cI}_{\om},$$
that is,
$\overline{\cI}_+=\cup_{\om\in\cT'(\overline{\cI}_+)}\overline{\cI}_{\om}$. Hence for any set $\Omega''$ such that $\cT'(\overline{\cI}_+)\subseteq\cT''\subseteq\cT_+$, we have $\overline{\cI}_+=\cup_{\om\in\cT''}\overline{\cI}_{\om}$. Let $\cT'(\overline{\cI})=\cT'(\overline{\cI}_+)\cup\cT'(\overline{\cI}_-)$. Then $\overline{\cI}=\overline{\cI}_+\cup\overline{\cI}_-=\cup_{\om\in\cT'(\overline{\cI})}\overline{\cI}_{\om}$. Similarly, $$\overline{\cO}_+=\underset{\om\in\cT'(\overline{\cO})_+}{\cup}\overline{\cO}_{\om},\quad
 \overline{\cO}_-=\underset{\om\in\cT'(\overline{\cO})_-}{\cup}\overline{\cO}_{\om},\quad
  \overline{\cO}=\underset{\om\in\cT'(\overline{\cO})}{\cup}\overline{\cO}_{\om},$$ for some finite $\cT'(\overline{\cO})\subseteq\cT$, where $\cT'(\overline{\cO})_{\pm}=\{\omega\in\cT'(\overline{\cO})\colon {\rm sgn}(\om_1)=\pm1\}$. By the definition of gcd, $\om^*$ is a common divisor of any subset of $\cT$, and there exists $-m_1\om^*\in\cT_-$ and $m_2\om^*\in\cT_+$ with $m_1,\ m_2\in\N$ coprime.  Let $\widetilde{\cT}=\cT'(\overline{\cI})\cup\cT'(\overline{\cO})\cup\{-m_1\om^*,m_2\om^*\}$. Then $\gcd(\widetilde{\cT})=\om^*$ and $$\underset{\om\in\widetilde{\cT}_+}{\cup}\overline{\cI}_{\om}=\overline{\cI}_+,\ \underset{\om\in\widetilde{\cT}_-}{\cup}\overline{\cI}_{\om}=\overline{\cI}_-,\ \underset{\om\in\widetilde{\cT}}{\cup}\overline{\cI}_{\om}=\overline{\cI},\ \underset{\om\in\widetilde{\cT}_+}{\cup}\overline{\cO}_{\om}=\overline{\cO}_+,\ \underset{\om\in\widetilde{\cT}_-}{\cup}\overline{\cO}_{\om}=\overline{\cO}_-,\ \underset{\om\in\widetilde{\cT}}{\cup}\overline{\cO}_{\om}=\overline{\cO},$$
 where $\widetilde{\cT}_{\pm}=\{\omega\in\widetilde{\cT}\colon {\rm sgn}(\om_1)=\pm1\}$.
\pre


\begin{thebibliography}{9}
\footnotesize
\bibitem{ABCJ18}
\textsc{Anderson, D.F.}, \textsc{Brijder, R.}, \textsc{Craciun, G.}, and \textsc{Johnston, M.D.} (2018).
Conditions for extinction events in chemical reaction networks with discrete state spaces.
\textit{J. Math. Biol.}
\textbf{76}, 1535--1558.


\bibitem{ACKK18}
\textsc{Anderson, D.F.}, \textsc{Cappelletti, D.}, \textsc{Koyama, M.}, and \textsc{Kurtz, T.G.} (2018).
Non-explosivity of stochastically modeled reaction networks that are complex balanced.
\textit{Bull. Math. Biol.}
\textbf{80}, 2561--2579.


\bibitem{AEJ14}
\textsc{Anderson, D. F.}, \textsc{Enciso, G.}, and \textsc{Johnston, M.D.} (2014).
Stochastic analysis of biochemical reaction networks with absolute concentration robustness.
\textit{J. R. Soc. Interface}
\textbf{11}, 20130943.

\bibitem{AK11}
\textsc{Anderson, D.F.} and  \textsc{Kurtz, T.G.} (2011).
Continuous time Markov chain models for chemical reaction networks.
In: \textsc{Koeppl, H.}, \textsc{Setti, G.}, \textsc{di Bernardo, M.}, and \textsc{Densmore, D.} (eds),
\textit{\em Design and Analysis of Biomolecular Circuits: Engineering Approaches to Systems and Synthetic Biology}.
Springer-Verlag, New York, 3--42.


\bibitem{C97}
\textsc{Chen, R.-R.} (1997).
An extended class of time-continuous branching processes.
\textit{J. Appl. Prob.}
\textbf{34}, 14--23.



\bibitem{CMM13}
\textsc{Collet, P.}, \textsc{Mart\'{i}nez, S.}, and \textsc{San Mart\'{i}n, J.} (2013).
\textit{Quasi-Stationary Distributions: Markov Chains, Diffusions and Dynamical Systems}.
Probability and its Applications (New York), Springer, Heidelberg.


\bibitem{E79}
\textsc{Ewens, W.J.} (2004).
\textit{Mathematical Population Genetics. I. Theoretical introduction}.
Interdisciplinary Applied Mathematics, Springer-Verlag. New York.

\bibitem{F19}
\textsc{Feinberg, M.} (2019).
\textit{Foundations of Chemical Reaction Network Theory}.
Applied Mathematical Sciences, Springer International Publishing. Cham.


\bibitem{G83}
\textsc{Gardiner, C.W.} (2009).
\textit{Stochastic Methods: A Handbook for Physics, Chemistry and the Natural Sciences}.
4th ed., Springer Series in Synergetics. Springer-Verlag, Berlin.

\bibitem{GH98}
\textsc{Gross, D.} and \textsc{Harris, C.M.} (1998).
\textit{Fundamentals of Queuing Theory}.
3rd ed., Wiley Series in Probability and Statistics, John Wiley $\&$ Sons, Inc., New York.

\bibitem{HW20}
\textsc{Hansen, M.C.}  and \textsc{Wiuf, C.} (2020).
Existence of a unique quasi-stationary distribution for stochastic reaction networks.
\textit{Electron. J. Probab.}
\textbf{25}, 1--30.





\bibitem{HJ72}
\textsc{Horn, F.} and \textsc{Jackson, R.} (1972).
General mass action kinetics.
\textit{Arch. Rational Mech. Anal.}
\textbf{47}, 81--116.

\bibitem{Klebaner}
\textsc{Klebaner, F.C.},  and \textsc{Liptser, R.} (2001).
Asymptotic analysis and extinction in a stochastic Lotka-Volterra model.
\textit{Ann. Appl. Probab.}
\textbf{11}, 1263--1291.


\bibitem{N98}
\textsc{Norris, J.R.} (1998).
\textit{Markov Chains}.
Cambridge Series in Statistical and Probabilistic Mathematics, \textbf{2}, Cambridge Univ. Press, Cambridge, UK.


\bibitem{PCMV15}
\textsc{Pastor-Satorras, R.}, \textsc{Castellano, C.}, \textsc{Van Mieghem, P.}, and \textsc{Vespignani, A.} (2015).
Epidemic processes in complex networks.
\textit{Rev. Mod. Phys.}
\textbf{87}, 925--979.


\bibitem{PCK14}
\textsc{Paulev\'{e}, L.}, \textsc{Craciun, G.}, and \textsc{Koeppl, H.} (2014).
Dynamical properties of discrete reaction networks.
\textit{J. Math. Biol.}
\textbf{69}, 55--72.


\bibitem{R57}
\textsc{Reuter, G.E.H.} (1957).
Denumerable Markov processes and the associated contraction semigroups on $l$.
\textit{Acta Math.}
\textbf{97}, 1--46.

\bibitem{SF10}
\textsc{Shinar, G.} and \textsc{Feinberg, M.} (2010).
Structural sources of robustness in biochemical reaction networks.
\textit{Science}
\textbf{327}, 1389--1391.


\bibitem{WH83}
\textsc{Weidlich, W.} and \textsc{Haag, G.} (1983).
\textit{Concepts and Models of a Quantitative Sociology: The Dynamics of Interacting Populations}.
Springer Series in Synergetics, Springer-Verlag, Berlin.



\bibitem{W06}
\textsc{Wilkinson, D.J.} (2006).
\textit{Stochastic Modelling for Systems Biology}.
Chapman $\&$ Hall/CRC Mathematical and Computational Biology Series, Chapman $\&$ Hall/CRC, London, UK.


\bibitem{WX20}
\textsc{Wiuf, C.} and \textsc{Xu, C.} (2020).
Classification and threshold dynamics of stochastic reaction networks.
\textit{Preprint}. arXiv:2012.07954.

\bibitem{XHW20a}
\textsc{Xu, C.}, \textsc{Hansen, M.C.}, and \textsc{Wiuf, C.} (2020).
Full classification of dynamics for one-dimensional continuous time Markov chains with polynomial transition rates.
\textit{Preprint}. arXiv:2006.10548.

\bibitem{XHW20b}
\textsc{Xu, C.}, \textsc{Hansen, M.C.}, and \textsc{Wiuf, C.} (2020).
The asymptotic tails of limit distributions of continuous time Markov chains.
\textit{Preprint}. arXiv:2007.11390.

\end{thebibliography}
\end{document}